\documentclass{amsart}
\input xypic
\usepackage{amssymb}   
\usepackage{amsmath}

\newtheorem{Theorem}{Theorem}[section]
\newtheorem{Lemma}[Theorem]{Lemma}
\newtheorem{Proposition}[Theorem]{Proposition}

\theoremstyle{definition}
\newtheorem{Definition}{Definition}
\newtheorem{Remark}[Theorem]{Remark}

\newtheorem{example}{Example}

\newcommand{\Sing}{\operatorname{Sing}}

\newcommand{\mult}{\operatorname{mult}}

\renewcommand{\O}{{\mathcal O}}

\newcommand{\I}{{\mathcal I}}
\newcommand{\E}{{\mathcal E}}

\newcommand{\p}{{\mathbb P}}

\newcommand{\codim}{\operatorname{codim}}
\newcommand{\red}{\operatorname{red}}

\newcommand{\Sec}{\operatorname{Sec}}

\newcommand{\rk}{\operatorname{rk}}

\newcommand{\map}{\dasharrow}

\begin{document}

\title[Some extremal contractions]{Some  extremal contractions between smooth
varieties arising from projective geometry}
\author{Alberto Alzati}\thanks{First author partially supported by the
M.I.U.R. of the Italian Government
 in the framework of the national Research Project ``Geometry of algebraic varieties"
 (Cofin 2002).}
\address{Dipartimento di Matematica ``F. Enriques",
 Universit\` a degli
Studi di Milano, Via C. Saldini 50, 20133 Milano, Italia}
\email{alzati@mat.unimi.it}

\author{Francesco Russo}
\address{Departamento de Matematica, Universidade Federal
de Pernambuco, Cidade Universitaria, 50670-901 Recife--PE, Brasil}
\email{frusso@dmat.ufpe.br}
\thanks{Second author  partially supported by CNPq
(Centro Nacional de Pesquisa), grant 300961/2003-0 and by
PRONEX-Algebra Comutativa e Geometria Algebrica}

\maketitle

\section*{Introduction}
Let $\psi:X\to Y$ be a proper morphism with connected fibers from
a smooth projective variety $X$ onto a normal variety $Y$,
 i.e.  a {\it contraction}. If $-K_X$ is
$\psi$-ample, then $\psi$ is said to be an {\it extremal
contraction} and if moreover $Pic(X)/\psi^*(Pic(Y))\simeq
\mathbb{Z}$ then $\psi$ is said to be an {\it elementary extremal
contraction} or a {\it Fano-Mori contraction}. Usually
contractions are divided into two types: birational and of {\it
fiber type}, i.e. with $\dim(X)>\dim(Y)$.

The importance of this study and the parallel development of the
notion of extremal ray and of the possible contractions of these
rays was established by Mori in
 \cite{Mori}.
In dimension 2 and over an arbitrary field only three types of
Fano-Mori contractions exist: the inverse of the blow-up of a
reduced and geometrically irreducible set of  points defined over
the field (birational), conic bundle structure and del Pezzo
surfaces (of fiber type) and this classification is essentially
equivalent to the theory of minimal models in dimension 2. In
dimension 3 the theory is more delicate but there is a complete
description for algebraically closed fields given by Mori in
characteristic zero and by Koll\' ar for arbitrary algebraically
closed fields, see for example \cite{KM}, theorem 1.32 or
\cite{Mori}. The common feature of extremal elementary
contractions in dimension 3 (and 2) is that if $E\subseteq X$ is
the exceptional locus of the contraction, i.e. the locus of
positive dimensional fibers, then $\phi_{|E}:E\to\phi(E)\subseteq
Y$ is equidimensional. Moreover in the birational case (and over
an algebraically closed field) these contractions are always the
inverse of the blow-up of a (possible singular) point or of a
smooth curve, see {\it loc. cit.}

In dimension greater than three the situation is more complicated
as it is shown by simple examples and no general result is known.
Here we furnish some examples of Fano-Mori contractions $\psi:X\to
Y$ between smooth projective varieties with $\dim(X)\geq 4$
defined over "sufficiently" small fields and in arbitrary
characteristic, having  some exceptional fibers of unexpected
dimension. All these contractions come from classical projective
geometry, or more precisely they are associated to the projective
geometry of remarkable classes of varieties: varieties with one
apparent double, triple, quadruple point or varieties defining
{\it special} Cremona transformations. Surely many examples of the
phenomenon of non-equidimensionality of the fibers of $\phi_{|E}$
were constructed  and are well known at least in characteristic
zero or over algebraically closed fields (see for example
\cite{CKM}, pg. 35, or \cite{KM}, pg. 44, and also \cite{AW1} or
\cite{AW2}). In the examples presented below the accent is on the
geometric nature of these constructions and on the close relations
with classical results. Moreover, they are also the first
instances, from a historical point of view, of a class of
morphisms with isolated exceptional fibers, a phenomenon studied
in detail, in a  local setting, in \cite{AW1}, \cite{AW2}. In
these last two papers some examples of these phenomena were
constructed over the complex field by using vector bundles and
geometrical constructions.

In section 2 we  construct interesting Fano-Mori contractions by
considering codimension 2 determinantal smooth subvarieties of
$\p^m$, $3\leq m\leq 5$, see proposition \ref{determinantal}.
Then we  furnish an example  of a Fano-Mori birational contraction
between smooth varieties not having equidimensional exceptional
fibers (so that it is not the inverse of a blow-up with smooth
center) by considering  a very particular  special quarto-quartic
Cremona transformation of $\p^4$ studied classically in the
general case by Todd and Room, see \cite{Todd} and \cite{Room}.
The existence of this exceptional quarto-quartic special Cremona
transformation seemed not to be noticed before, see example
\ref{surfp4}.  We remark that the general quarto-quartic
transformation does not give this exceptional phenomenon. The
particular  contraction is related to the construction of a
determinantal degree 10 smooth surface containing plane quartic
curves, whose existence was essentially suggested by Room, see
example $iii)$ on page 80 of \cite{Room}.

An example of not  equidimensional fiber
type Fano-Mori contraction can be constructed from the congruence
of trisecant lines to a Bordiga surface in $\p^4$, a surface with
one apparent triple point according to Severi (\cite{Sev}). The
double-ten of points-planes associated to it (see page 71 of
\cite{Room}) shows that the associated contraction, which is of
fiber type in dimension 4, has 10 exceptional two dimensional
fibers, the focal locus of the congruence outside the Bordiga
surface, see example \ref{Bordiga}. This example is in a certain
sense also more classical and its origins are harder to trace.

In dimension 5 the determinantal construction furnishes an example of
fiber type Fano-Mori contraction onto $\p^4$ coming from a degree
10 determinantal 3-fold $Z\subset \p^5$ having one apparent
quadruple  point. In the general case the congruence of its
4-secant lines has  a smooth curve of exceptional 2 dimensional
fibers. In a particular constructed example it has a curve  of
exceptional 2 dimensional fibers and exactly one 3 dimensional
fiber corresponding to the unique  singular point of the curve,
see example \ref{conicp5}. A general hyperplane section of this
particular example furnishes a family  of surfaces (and hence of
interesting Fano-Mori contraction) of the type discussed in
example \ref{surfp4}.

From some very particular special quinto-quintic Cremona
transformation of $\p^5$ we construct some Fano-Mori birational
contractions between smooth 5-folds and $\p^5$ having general
exceptional fiber of dimension 1 and a finite number of isolated
exceptional fibers of dimension 2 or 3, example \ref{birp5}.

Once the path has been cleared we continued the study and
furnished examples of different nature, whose interest appears
also in other important problems. Continuing the study in
dimension 4 we showed how the analysis given by Fano in
\cite{Fano1} of general smooth cubic hypersurfaces in $\p^5$
passing through a del Pezzo quintic surface can be extended to
arbitrary (infinite) fields and also to every smooth cubic
hypersurface through a del Pezzo surface. Some particular smooth
cubic hypersurfaces generate Fano-Mori birational contractions
between these smooth 4-folds and $\p^4$ with general exceptional
fiber of dimension 1 and with 1,2,3 or 4 exceptional fibers
isomorphic to $\p^2$. These fibers contract to singular points of
the fundamental locus of the inverse map, which is the projection
of a degree 9 normal (K3) surface of genus 8, see \ref{th:Fano}.
This result and the analogous for cubic hypersurfaces through a
degree 4 rational normal scroll in $\p^5$ are probably well known
and were rediscovered many times. We take the opportunity of
giving the exact classical reference and the extension to
arbitrary fields.

In dimension 6 an interesting  Fano-Mori birational contraction
can be constructed from Semple and Tyrrell example of a
quadro-quartic special Cremona transformation in $\p^6$
(\cite{ST2}), given by quadrics through a general rational octic
surface $S\subset\p^6$. In this case the associated  Fano-Mori
contraction $\psi:Bl_S(\p^6)\to \p^6$ is very interesting and has
a rather intricate behaviour for the exceptional fibers of
dimension 2: there are 28 isolated {\it planes} and a two
dimensional family of non-isolated fibers which are {\it
quadrics}, parametrised by a smooth quadric surface, see
\ref{SempleTyrrell}. Therefore there are isolated two dimensional
fibers and a two dimensional family of exceptional two dimensional
fibers in the same Fano-Mori contraction and therefore the local
behaviour and the global behaviour of exceptional fibers can be
completely different for the same elementary extremal contraction.
Probably this is the most intriguing example of Fano-Mori
birational contraction together with the next examples related to
it. Indeed, one constructs a
 3-fold $X\subset\p^7$ of degree 8 having  one apparent double point
 and having $S$ as its general hyperplane section;  this 3-fold produces
a Fano-Mori contraction of fiber type $\psi:Bl_X(\p^7)\to\p^6$ for
which the behaviour of the exceptional fibers of dimension greater
or equal to $2$ is very interesting and it is  treated completely
by geometrical methods, see \ref{Semplep7}. We conclude by
constructing a Fano-Mori birational contraction between smooth
cubic hypersurfaces  $Y\subset\p^7$ through the above mentioned
3-fold $X$ and $\p^6$ and we study in detail the inverse map and
all the exceptional fibers. In particular we show that
$\phi^{-1}:\p^6\map Y$ is given by quintic forms vanishing along a
4-fold $Z\subset\p^6$ of degree 13 and having at least a two
dimensional locus of double points, over which the fiber of the
contraction $\widetilde{\phi}:Bl_X(Y)\to\p^6$ consists of
(non-isolated) quadric surfaces, see \ref{Fanop7}.

This list can be considered as a natural continuation of the
classical study of examples of interesting geometrical phenomena.
The guiding principle was the beginning of Semple and
Tyrrell paper \cite{ST2}: {\it Any Cremona transformation defined
by a homaloidal system of primals with a single irreducible
non-singular base variety is a rare enough phenomenon to merit
special study.}
 We were deeply influenced by the classical
papers and books \cite{Sev}, \cite{ST2}, \cite{ST1}, \cite{ST3},
\cite{SR}, \cite{Room}, \cite{Todd}, \cite{Fano1} and the more
recent \cite{ESB}.

We conclude by remarking that even if these examples could have
some interest for higher dimensional geometry from the point of
view of contractions, the methods we used are completely
elementary and geometric and all the results can be considered
classical or generalizations of them. Last but not least, we claim
no originality for the geometrical interpretations of these
classical examples.
\medskip

{\bf Acknowledgements}. The second author would like to thank
Professor Shepherd-Barron for asking him about the classical
study of the general contractions in  \ref{th:Fano}. This was the
starting point of this list of examples.

\section{Definitions and  preliminary results}
Let us recall the following definitions.
\begin{Definition}
{\rm Let $\phi:\p^n\map\p^n$ be a birational map, which is  not an
isomorphism. Then the linear system $\phi^{-1}_*|\O(1)|$ is called
a homaloidal system. If the map $\phi$ is given by the homogeneous
forms $F_0,\ldots,F_n$ of degree $d\geq 2$ defining a smooth,
geometrically irreducible variety
$X^1=V(F_0,\ldots,F_n)\subset\p^n$, then $\phi$ is said a {\it
special Cremona transformation} and the associated linear system a
{\it special homaloidal system}. More generally, a linear system
of hypersurfaces of $\p^r$ is said to be {\it subhomaloidal} if
the (closure of a) general fiber of the associated rational map
$\phi:\p^r\map\p^n$, $r\geq n$, is a linear projective space. It
is said to be {\it special} if the base locus scheme of the linear
system is  smooth and irreducible;
 it is said to be {\it completely subhomaloidal} if (the
closure of) every fiber is a linear projective space. These
definitions can also be extended to linear systems of
hypersurfaces such that the associated map $\phi:\p^r\map
Z\subseteq\p^m$, $Z=\overline{Im(\phi)}$, is birational onto the
image, respectively has a linear projective space as general fiber. }

\end{Definition}

\medskip

Thus, let $\phi:\p^n\dasharrow \p^n$ be a special Cremona
transformation  of type $(d_1,d_2)$, which means that $\phi$ is
given by a linear subsystem of $|\O(d_1)|$ and $\phi^{-1}$ by a
linear subsystem of $|\O(d_2)|$. Let $X^1$ be the base locus
scheme of $\phi$, which is smooth and geometrically irreducible by
definition, and let $X^2$ be the base locus scheme of $\phi^{-1}$.
Let $r_i=\dim(X^i)$ for $i=1,2$. Let $\pi:Bl_{X^1}(\p^n)\to \p^n$
stand for the structural map. We can consider $Bl_{X^1}(\p^n)$
contained in $\p^n\times\p^n$ so that there is a canonical diagram
$$\xymatrix{
          Bl_{X^1}(\p^n)\ar[d]_{\pi}\ar[dr]^{\psi}       &             \\
  \p^n  \ar@{-->}[r]_{\phi} & \p^n              }
$$

where $\psi$ is the projection onto the second factor.

Put $E_1=\pi^{-1}(X^1)$ and $E_2=\psi^{-1}(X^2)$. Take $H_1\in
|\pi^*(\O(1))|$, $H_2\in |\psi^*(\O(1))|$ and let $T_1$ and $T_2$
stand, respectively, for the strict transform on
$Bl_{X^1}(\p^n)=Bl_{X^2}(\p^n)$ of a general line contained in
$\p^n$ and of a general line of the target  $\p^n$.

Let us recall that, given a subvariety $W\subset \p^n$ and an
integer $s\geq 2$, an {\it $s$-secant line} to $W$ is a line in
$\p^n$ not contained in $W$ and intersecting $W$ in at least $s$
points (counted with multiplicities over $\overline{k}$). The
closure $\mbox{\rm Sec}_s(W)\subset \p^n$ of the union of all
$s$-secant lines to $W$ is the {\it variety of $s$-secant lines}.
Of course, $\mbox{\rm Sec}_2(W)$ is the ordinary secant variety of
$W$.

The following lemma contains the preliminary  results of
\cite{ESB} and it will be a fundamental tool for our investigation
($\sim$ stands for linear equivalence).

\begin{Lemma}{\rm (\cite{ESB}, section 1)}\label{equations}
Let notations be as above.
\begin{enumerate}
 \item[{\rm (i)}] $H_2\sim d_1H_1-E_1$,
$H_1\sim d_2H_2-E_2$,  $E_2=(E_2)_{\red}$ and $T_1\cdot H_2=d_1$,
$T_2\cdot H_1=d_2$, $E_1\cdot T_2=E_2\cdot T_1=d_1d_2-1$.
Moreover, $E_2\sim (d_1d_2-1)H_1-d_2E_1$.

\item[{\rm (ii)}] $\psi(E_1)=\mbox{\rm Sec}_{d_2}(X^2)$ is a
hypersurface of degree $d_1d_2-1$ and $\pi(E_2)=\mbox{\rm
Sec}_{d_1}(X^1)$ is a hypersurface of degree $d_1d_2-1$.
\item[{\rm (iii)}] $2+r_1=d_2[(n+1)-d_1(n-r_1-1)]$.

\item[{\rm (iv)}] $2+r_2=d_1[(n+1)-d_2(n-r_2-1)]$.

\item[{\rm (v)}] $h^0(\I_{X^1}(d_1))=n+1$,
$h^0(\I_{X^1}(d_1-1))=0$, $h^0(\I_{X^1}(1))=0$ and $X^1$ is not a
complete intersection.

\item[{\rm (vi)}] If the base field $k$ is algebraically closed,
then $\psi^{-1}(U)\to U$  is the blow-up of $U:=\p^n\setminus
\Sing((X^2)_{\red})$ along $(X^2)_{\red}\cap U$. In $U$,
$(X^2)_{\red}=X^2$ so that $X^2$ is generically smooth.
\end{enumerate}
\end{Lemma}
\medskip

We recall the following theorem of Danilov.

\begin{Theorem}{\rm (\cite{Dan})} Let $f:X\to Y$ be a projective birational
morphism of smooth
algebraic spaces defined over an arbitrary field $k$  such that
$\dim(f^{-1}(y))\leq 1$ for every $y\in Y$. Then $f$ can be
decomposed into a sequence of monoidal transformations with smooth
centers of codimension 2.\label{te:Danilov}
\end{Theorem}
\medskip

The next remark collects  some general facts for further
reference.
\medskip

\begin{Remark} {\rm Let $X\subset\p^m$, $m\geq n+2$, be an
irreducible non-degenerate variety of dimension $n$. It is well
known that $n+2\leq \dim(\Sec_2(X))\leq\min\{m,2n+1\}$. We recall
that, if the base field has characteristic zero, the main result
of \cite{Ran} assures $\dim(\Sec_{n+2}(X))\leq n+1$. This last
property is true over an arbitrary field for smooth irreducible
curves by the Trisecant Lemma, see for example \cite{Ha} IV.3, and
also for the examples of varieties $X\subset\p^{n+2}$ of section 2
by a direct verification, using the following elementary
observation.

Suppose that for a general point $p\in \Sec_s(X)\subseteq\p^m$,
$s\geq 2$, the union of the $s$-secant lines through $p$ is a
projective space $\p^r_p\subseteq\p^m$, $r\geq 2$, cutting scheme
theoretically $X$ along a hypersurface in $\p^r_p$. Let
$X'\subset\p^{m-1}$ be a general hyperplane section of $X$. Then
for a general point $p\in \Sec_s(X')$, the union of the $s$-secant
lines through $p$ is a projective space
$\p^{r-1}_p\subseteq\p^{m-1}$, cutting scheme theoretically $X'$
along a hypersurface in $\p^{r-1}_p$. }\label{re:secant}

\end{Remark}
\medskip

We conclude the section by combining all the above result in the
following proposition.

\begin{Proposition} Let $\phi:\p^n\map\p^n$ be a special Cremona
transformation of type $(d_1,d_2)$ and let notations be as above.
Then $\psi:\widetilde{\p^n}=Bl_{X^1}(\p^n)\to\p^n$ is an extremal
elementary contraction between smooth varieties. The morphism
$\psi$ is the inverse of a blow-up with smooth center if and only
if $X^2$ is smooth, i.e. if and only if $\phi^{-1}$ is special. If
$\phi^{-1}$ is not special and if the base field is algebraically
closed, then $\psi:\psi^{-1}(U)\to U$ is the blow-up of
$U=\p^n\setminus Sing((X^2)_{\red})$ along $X^2\setminus
\Sing((X^2)_{\red})$. If $\phi^{-1}$ is not special and
$\codim(X^2)=2$, then $\psi:\psi^{-1}(V)\to V$ is the blow-up of
$V:=\p^n\setminus \{x\in X^2\;:\; \dim(\psi^{-1}(x))\geq 2\}$
along $X^2\setminus \{x\in X^2\;:\; \dim(\psi^{-1}(x))\geq
2\}$.\label{bircontraction}
\end{Proposition}
\begin{proof} By Zariski's Main theorem and by Zariski's
connectedness theorem, the morphism $\psi$ has connected fibers
and it is elementary by lemma \ref{equations}. Now we show that
$-K_{\widetilde{\p^n}}$ is $\psi$-ample. Since $\psi$ is given by
$|d_1H_1-E_1|$, on an arbitrary fiber $F$ of $\psi$ we have
$d_1H_{1|F}=E_{1|F}$. Then from lemma \ref{equations} part iv) we
get $-K_{\widetilde{\p^n}|F}=(n+1)H_{1|F}-(n-r_1-1)E_{1|F}=
(n+1-d_1(n-r_1-1))H_{1|F}=\frac{2+r_1}{d_2}H_{1|F}$. The last
parts follow from lemma \ref{equations} part vi), if
$k=\overline{k}$, or by theorem \ref{te:Danilov} in the remaining
case; recall that in the last case by the theorem of the dimension
of the fibers $X^2\setminus \{x\in X^2\;:\; \dim(\psi^{-1}(x))\geq
2\}$ is not empty.
\end{proof}

\section{Determinantal examples}

In this section we will construct some  examples of Fano-Mori
contractions both birational and of fiber type by considering
some codimension 2 determinantal smooth subvarieties in $\p^4$
and in $\p^5$ having an interesting and rich extrinsic geometry.

Let $A(\mathbf{x})$ be a $n+1\times n$ matrix of linear forms on
$\p^m$ with coordinates $\mathbf{x}=(x_0:\ldots:x_m)$, $m\geq n$.
Suppose the $n+1$ minors of maximal order $n$ of $A(\mathbf{x})$,
$F_0,\ldots,F_n$ define a smooth geometrically irreducible
codimension 2 subvariety $X_{m-2}\subset\p^m$. By definition, the
variety $X_{m-2}$ is scheme theoretically defined by the forms
$F_i$'s and, for example by using   the Hilbert-Burch theorem (or
the Eagon-Northcott complex), the ideal sheaf $\I_{X_{m-2}}$ has
the minimal resolution
$$0\to\O(-n-1)^n\stackrel{A(\mathbf{x})}\to\O(-n)^{n+1}\to\I_{X_{m-2}}\to
0.$$ In particular $X_{m-2}$ is an arithmetically Cohen-Macaulay
variety. The homogeneous forms $F_0,\ldots,F_n$  define a
rational map $\phi:\p^m\map\p^n$.

Let $\mathbf{y}=(y_0:\ldots:y_n)$ be coordinates on the target
$\p^n$. The closure of the graph of $\phi$, $\Gamma_{\phi}$, in
$\p^m\times\p^n$ is $Bl_{X_{m-2}}(\p^m)=\widetilde{\p^m}$, which
coincides outside the exceptional divisor $E$ with the scheme of
equation $$\sum
a_{i,j}(\mathbf{x})y_j=A(\mathbf{x})^t\cdot\mathbf{y}^t=\mathbf{0}.$$

In fact for $\mathbf{x}\in\p^m\setminus X_{m-2}$ by elementary
linear algebra we get: $$A(\mathbf{x})^t\cdot
\mathbf{y}^t=\mathbf{0}\Leftrightarrow
\mathbf{y}=(F_0(\mathbf{x}):\ldots:F_n(\mathbf{x}))\Leftrightarrow
(\mathbf{x},\mathbf{y})\in\Gamma_{\phi}.$$

 There
exists a $n\times (m+1)$ matrix of linear forms $B(\mathbf{y})$
such that $$A(\mathbf{x})^t\cdot \mathbf{y}^t=B(\mathbf{y})\cdot
\mathbf{x}^t.$$ Let $\widetilde{\phi}:\widetilde{\p^m}\to\p^n$ be
the resolution of the rational map $\phi$. Let
$\mathbf{y}=\phi(\mathbf{x})$, $\mathbf{x}\in \p^m\setminus
X_{m-2}$. The above discussion yields
$$\phi^{-1}(\mathbf{y})=\{\;\;\mathbf{x}\in\p^m\setminus X_{m-2}\; :\;
B(\mathbf{y})\cdot \mathbf{x}^t=\mathbf{0}\;\;\},$$ whose closure
in $\p^m$ is clearly a linear space of dimension
$m-\rk(B(\mathbf{y}))$, let us say
$\p_{\mathbf{y}}^{m-\rk(B(\mathbf{y}))}$. In particular we deduce
that the fibers of $\phi$ are all irreducible and that the rank of
$B(\mathbf{y})$ at the general point of $Im(\phi)$ is equal to
$\dim(Im(\phi))$ by the theorem of the dimension of the fibers.

It is also clear that, when
$\p_{\mathbf{y}}^{m-\rk(B(\mathbf{y}))}$ is positive dimensional,
it cuts $X_{m-2}$ along a determinantal hypersurface of degree
$n$. We leave the details to the reader. The geometrical
interpretation is that if through a point $p\in\p^m\setminus
X_{m-2}$ there passes a $n$-secant line to $X_{m-2}$, since this
line is contracted by $\phi$, the locus of $n$-secant lines to
$X_{m-2}$ through $p$ is the projective space
$\p_{\phi(p)}^{m-\rk(B(\phi(p))}$ cutting $X_{m-2}$ along a
hypersurface whose equation is determinantal in this projective
space. Moreover this geometric description implies that, by
looking at $\widetilde{\phi}$, we can say that every fiber of
$\widetilde{\phi}$ is either a linear space, eventually of varying
dimension, or it is isomorphic to a subvariety of a linear space
of the form $\p^{r(q)}\times q$, $q\in\p^n$,  contained in the
exceptional divisor $E$. We will show below that $\phi$ is always
dominant. If $m=n$, then $\phi$ is birational and it is a
Fano-Mori contraction by proposition \ref{bircontraction}. If
$m>n$  the same argument used in the proof of proposition
\ref{bircontraction} easily gives that $-K_{\widetilde{\p^m}}$ is
$\widetilde{\phi}$-ample. In any case,  $\widetilde{\phi}$ is a
Fano-Mori elementary contraction.

If $3\leq m\leq5$ it is known that examples of smooth varieties
$X_{m-2}$ as above exist in characteristic zero by Bertini's
theorem for spanned vector bundles. It is not difficult to
construct smooth, geometrically irreducible  examples over
suitable "small" fields, where the above construction applies as
we see below. If $m\geq 6$ the above arguments apply but the
variety $X_{m-2}$ will be singular. For $m=n$ and for
$A(\mathbf{x})$ sufficiently general (for example of rank $n$ and
with good Fitting ideals), one obtains a series of Cremona
transformations of type $(n,n)$ investigated by classical
geometers and in particular by Godeaux, see \cite{Godeaux} and
also \cite{RS}, example 2.4.

The above discussion about the locus of $n$-secant lines through a
general point of $\p^m$ is a particular version of proposition
2.9.1 of \cite{Room}. This result of Room can be considered as the
first case of the general situation described in proposition
\ref{prop:Vermeire}, see also \cite{HKS}, \cite{Ve} and \cite{RS}.
In fact Room's result is true for more general situations and
determinantal varieties, see {\it loc. cit.}

If $m=n$ (and if $3\leq n\leq 5$), it is easy to see that
$\Sec_n(X_{n-2})\subsetneqq\p^n$, for example by applying
\ref{re:secant} to $X_{n-2}$ or to its general hyperplane section.
Therefore $\phi$ will be dominant and birational, i.e.
$\phi:\p^n\map\p^n$ will be a special Cremona transformation of
type $(n,n)$, by lemma \ref{equations} iii), whose base locus is
exactly $X_{n-2}$. The general fiber is necessarily a point
because $\Sec_n(X_{n-2})\subsetneqq\p^n$, so that
$\rk(B(\mathbf{y}))=n$ for a general point $\mathbf{y}\in
Im(\phi)$. By applying the same argument to the matrix of linear
forms $B(\mathbf{y})$ and from the fact that $\phi^{-1}$ is given
by forms of degree $n$, we get that the inverse map is given by
the maximal minors of $B(\mathbf{y})$. The inverse map is not
defined along a determinantal scheme $X^2$ of the same type, which
in general will not be smooth but only generically smooth.  From
this description it is also clear that $\Sing(X^2)$ will be given
by the vanishing of the $(n-1)\times(n-1)$ minors of
$B(\mathbf{y})$ so that it is easily computable. One can look at
remark 2.3 of \cite{RS} for an example of non-special Cremona
transformation of similar type for which the maximal minors of
$B(\mathbf{y})$ have a common factor so that the corresponding
$\phi$ is not of type $(n,n)$.

The proper transform  of $\Sec_n(X_{n-2})$ on $\widetilde{\p^n}$
will be the geometrically irreducible divisor $E_2$ by lemma
\ref{equations}. The general fiber of
$\widetilde{\phi}_{|E_2}:E_2\to X^2$ has dimension 1 by the
theorem of the fibers and it is an $n$-secant line to $
X^1=X_{n-2}$. By lemma \ref{bircontraction},
$\widetilde{\phi}:\widetilde{\phi}^{-1}(V)\to V$ is the blow-up of
$V=\p^n\setminus \{x\in X^2\;:\; \dim(\psi^{-1}(x))\geq 2\}$ along
$X^2\setminus \{x\in X^2\;:\; \dim(\psi^{-1}(x))\geq 2\}$. By the
above description in this case $\{x\in X^2\;:\;
\dim(\psi^{-1}(x))\geq 2\}=\Sing(X^2)$.

If $m=n+1$ (and to have smooth $X_{n-1}=X^1\subset\p^{n+1}$,
$2\leq n\leq 4$),  the description of the fibers of $\phi$ yields
$\Sec_n(X^{1})=\p^{n+1}$ and the fact that the general fiber $F$
of the associated rational map $\phi:\p^{n+1}\map\p^n$ is a
$n$-secant line, see remark \ref{re:secant}. The contraction
$\widetilde{\phi}$ is easily seen to be extremal and elementary by
arguing as in proposition \ref{bircontraction}. Moreover, if $F$
is a one dimensional fiber as above,
$-K_{\widetilde{\p^m}_{|F}}\simeq \O_{\p^1}(2)$. With the
classical terminology, the variety $X^1$ {\it has one apparent
n-ple point}. Take now an arbitrary hypersurface
$Y=V(G)\subset\p^{n+1}$ of degree $n+1$ through
$X_{n-2}\subset\p^{n+1}$. The restriction of $\phi$ is a
birational isomorphism $\phi_{|Y}:Y\map\p^{n}$ because a general
fiber of $\phi$ is a $n$-secant line to $X_{n-2}$ not contained in
$Y$, so that this general fiber residually cuts $Y$ in a point
outside $X_{n-2}$. Moreover the equation of an $Y$ as above can be
written as the determinant of a $(n+1)\times(n+1)$ matrix of
linear forms on $\p^{n+1}$. In fact since
$G(\mathbf{x})=\sum_{i=0}^nl_i(\mathbf{x})F_i(\mathbf{x})$, the
graph of $\phi_{|Y}$ inside $\p^{n+1}\times\p^n$ coincides outside
$E\cap Y$ with the scheme given by the $n+1$ bihomogeneous $(1,1)$
forms
$$A(\mathbf{x})^t\cdot \mathbf{y}^t=\mathbf{0},$$
$$\sum_{i=0}^nl_i(\mathbf{x})y_i=0, $$
which can be written as
\begin{equation}\widetilde{A}(\mathbf{x})^t\cdot
\mathbf{y}^t=\mathbf{0},\end{equation} where
$\widetilde{A}(\mathbf{x})$ is the desired $(n+1)\times(n+1)$
matrix of linear forms on $\p^{n+1}$. The conclusion easily
follows since $\det(\widetilde{A}(\mathbf{x}))$ is a homogeneous
polynomial of degree $n+1$ vanishing at the general point of $Y$.
We remark that the same construction applies also for a
geometrically reduced, not necessarily geometrically irreducible,
determinantal variety $X_{n-1}\subset\p^{n+1}$, as soon as
$\Sec(X_{n-1})=\p^{n+1}$. By reverting the construction we have
that a hypersurface of degree $m$ in $\p^{m}$ containing a
codimension 2  geometrically reduced variety of the above type has
an equation given by the determinant of a $m\times m$ matrix of
linear forms on $\p^m$.

We collect the above discussion in the following proposition.
\begin{Proposition} Let $A(\mathbf{x})$ be a $n+1\times n$ matrix
of linear forms on $\p^m_k$, $k$ an arbitrary filed, with
coordinates $\mathbf{x}=(x_0:\ldots:x_m)$, $m\geq n$. Suppose the
$n+1$ minors of maximal order $n$ of $A(\mathbf{x})$,
$F_0,\ldots,F_n$ define a smooth geometrically irreducible
codimension 2 subvariety $X_{m-2}\subset\p^m$. The homogeneous
forms $F_0,\ldots,F_n$ define a rational map $\phi:\p^m\map\p^n$,
which resolves to a Fano-Mori contraction
$\widetilde{\phi}:\widetilde{\p^m}:=Bl_{X_{m-2}}\p^m\to\p^n$.

If $3\leq m=n\leq 5$, then $\phi:\p^n\map\p^n$ is a special
Cremona transformation of type $(n,n)$ and
$\widetilde{\phi}_{\widetilde{\phi}^{-1}(V)}:\widetilde{\phi}^{-1}(V)\to
V$ is the blow-up of $V=\p^n\setminus \Sing(X^2)$ along
$X^2\setminus \Sing(X^2)$.

If $3\leq m\leq 5$, $n=m-1$, then
$\widetilde{\phi}:Bl_{X_{m-2}}\p^{m}\to\p^{m-1}$ is generically a
conic-bundle, i.e. a one dimensional fiber $F$ of
$\widetilde{\phi}$ is an $m-1$-secant line to $X_{m-2}$ and
$-K_{\widetilde{\p^m}_{|F}}\simeq \O_{\p^1}(2)$. Moreover any
hypersurface $Y\subset\p^m_k$ of degree $m$ through $X_{m-2}$ has
an equation given by the determinant of a $m\times m$ matrix of
linear forms on $\p^m_k$ (and the same is true for an arbitrary
determinantal $X_{m-2}\subset\p^m$) and the restriction of $\phi$
to $Y$ is birational.\label{determinantal}
\end{Proposition}
\medskip
\begin{Remark} We remark that the result concerning the rationality over
the base field and the determinantal nature of arbitrary
hypersurfaces through some  determinantal reduced codimension 2
varieties, whose proof is elementary and geometric, can be
considered as a suitable generalization of the following fact.
Over a field $k$, a cubic surface
 $S\subset\p^3_k$ is rational over $k$ and it has an equation given by the determinant of
 a $3\times 3$ matrix of linear forms if it contains a
 degree 3 non-degenerate reduced connected curve defined over $k$, which is easily
 seen to have ideal generated by minors of a $2\times 3$ matrix of linear forms
 on $\p^3_k$. Usually
 over an algebraically closed field one takes as the above curve
 either a twisted cubic or the union of two skew lines and a line
 intersecting them or the union of a line and a conic intersecting
 in one point or three concurrent lines.
Since  a smooth cubic surface over an algebraically closed field
is easily seen to contain  a degree 3 curve of the first three
kinds, the above proposition applies to this class of
hypersurfaces, a classical result known to Grassmann and Clebsch,
see \cite{Gr}, \cite{Cl}, \cite{SR}, pg. 123 and also \cite{Be},
6.4. The extension to arbitrary fields was considered in \cite{BS}
and also in \cite{Be}, 6.5. Surely for smooth cubic surfaces the
point is to show geometrically the existence of the curve without
knowing a priori the plane representation of the cubic surface.
The morphism to $\p^2$, i.e. the plane representation for the ones
which are not cones, is obtained by solving the linear system (1).
\end{Remark}
\medskip

We are now ready to begin the construction  of examples of
Fano-Mori contractions with exceptional behaviour. We begin by
discussing the following  special completely subhomaloidal linear
system only as an elementary application of the above methods and
also because we need this result in the last sections. As we see
below, all the other examples presented later will be only
suitable generalizations of this case.

\medskip

\begin{example} {\rm  {\it Quadrics through the Segre 3-fold in $\p^5$
or its general linear section and associated elementary
contractions}. Let $n=2$ and $m=5$ and let notations be as above.
Let $\phi:\p^5\map\p^2$ be the associated completely subhomaloidal
linear system, by abusing language, of quadrics through
$X^1=\p^1\times\p^2$. Modulo projective transformations we can
suppose that $A(\mathbf{x})$ is the following $3\times 2$ matrix
of linear forms on $\p^5$, whose minors give rise to $X^1$:

$$\left(
\begin{array}{cc}
x_0 &x_3\\
x_1 &x_4\\
x_2 &x_5

\end{array}
\right).$$

The $2\times 6$ matrix $B(\mathbf{y})$ is $$\left(
\begin{array}{cccccc}
y_0 &y_1&y_2&0&0&0\\
0&0&0&y_0&y_1 &y_2

\end{array}
\right),$$

so that every fiber of $\widetilde{\phi}:\widetilde{\p^5}\to\p^2$
is a $\p^3$ cutting $X^1$ along a  quadric surface.}\label{segre}
\end{example}
\medskip

\begin{example} {\rm (Todd--Room contraction). {\it A Fano-Mori  birational contraction
between a smooth projective 4-fold and $\p^4$ having an isolated
exceptional fiber of dimension 2}. Let $m=n=4$ in the above
construction with $A(\mathbf{x})\in M_{5\times 4}(k[\mathbf{x}])$,
$k$ to be chosen during the construction. If it exists,
$X^1\subset\p^4$ is a smooth degree 10 and genus 11 surface
defined over $k$, whose ideal is generated by 5 quartic forms
having the first syzygies generated by the linear ones and
defining a special Cremona transformation. A general
$A(\mathbf{x})$ will provide such an example but we aim at
discussing  all the cases defining a smooth $X^1$.

Let $\widetilde{\phi}:\widetilde{\p^4}=Bl_{X_1}(\p^4)\to\p^4$ be
the resolution of the associated special Cremona transformation
$\phi:\p^4\map\p^4$. By proposition \ref{determinantal} $\phi$ is
of type $(4,4)$ and  $\widetilde{\phi}^{-1}(V)\to V$ is the
blow-up of $V=\p^4\setminus \Sing(X^2)$ along $X^2\setminus
\Sing(X^2)$. Recall that $\Sing(X^2)=\{y\in\p^4\,:\,
\rk(B(\mathbf{y}))\leq 2\}.$ If  $\mathbf{y}\in \Sing(X^2)$, then
$\widetilde{\phi}^{-1}(\mathbf{y})$ is a plane in $\p^4\times
\mathbf{y}$ cutting $X^1$ scheme theoretically along a plane
quartic determinantal curve.
 The singular points of the base-locus scheme
$X^2$ are then in bijective correspondence with quartic plane
curves contained in $X^1$. This number is necessarily finite as it
is easy to see, eventually passing to $\overline{k}$.

Now we  show that the singular points of $X^2$ are of multiplicity
2 for a general quartic hypersurface containing $X^2$ and of
multiplicity 3 for $X^2$. In fact, the image of a general plane
passing through a singular point $p_0\in \Sing(X^2)$ under the
rational map induced by the linear system of quartics through
$X^2$, restricted to the plane, is linked $(4,4)$ to the union of
$X^1$ and the corresponding plane over $p_0$, so that it has
degree 5. It follows that the base locus of the plane system of
quartics consists of 16-5=11 points: $p_0$ with multiplicity $h$
and other $t$ points with multiplicity 1 such that $11=h^2 +t$. As
the dimension of that system of plane quartics is 5, it is easy to
see that $h=2$ and $t=7$ so that $p_0$ is a point of multiplicity
$\deg(X^2)-t$=10-7=3 for $X^2$.

The above analysis  also shows that the exceptional planes will be
points of multiplicity 2 for $E_2$, the blow-up of $\Sec_4(X^1)$
along $X^1$, and that the tangent cone at a singular point $p\in
\Sing(X^2)$ is a cone over a twisted cubic, in fact it has
dimension 2 and degree 3. This is in accordance with the
main-theorem of \cite{AW1}, pg. 256. Note that the tangent cone is
linked $(2,2)$ to the tangent plane at some point of a smooth
Bordiga surface passing through $p$ and linked $(4,4)$ to $X^2$.

By proposition  \ref{determinantal},
$\widetilde{\phi}:\widetilde{\p^4}\to\p^4$ is a Fano-Mori
contraction, which is  the inverse of a blow-up with smooth center
if and only if $\Sing(X^2)=\emptyset$.

We proceed by constructing an example of  a smooth surface $X^1$
such that $\Sing(X^2)=(0:0:0:0:1)$ and defined over "sufficiently"
large (or small, depending on the point of view) fields. To this
aim it is sufficient to construct $X^1$ as given by the matrix of
linear forms $A(\mathbf{x})$ \vskip 0.5cm

 $$\left(
\begin{array}{cccc}
-2x_{1}+x_{0} & -2x_{2}+x_{0} & 2x_{0} & -x_{1}-x_{4} \\
x_{3}+x_{0} & -x_{1}+x_{2} & x_{1}-2x_{3} & 2x_{2}-x_{3} \\
-x_{1}-x_{3} & -x_{3}+2x_{4} & -2x_{4}+x_{0} & x_{1}-x_{0} \\
-2x_{1}+x_{4} & -x_{2}-x_{0} & x_{2}+x_{3} & x_{2}+x_{4} \\
x_{3} & x_{4} & x_{3} & x_{4}
\end{array}
\right).$$ \vskip 0.5cm

 One verifies that $X^1$ is smooth and that
the matrix $B(\mathbf{y})$ is \vskip 0.5cm

$$\left(
\begin{array}{ccccc}
-2y_{0}-y_{2}-2y_{3} & 0 & y_{1}-y_{2}+y_{4} & y_{3} & y_{0}+y_{1} \\
-y_{1} & -2y_{0}+y_{1}-y_{3} & -y_{2} & 2y_{2}+y_{4} & y_{0}-y_{3} \\
y_{1} & y_{3} & -2y_{1}+y_{3}+y_{4} & -2y_{2} & 2y_{0}+y_{2} \\
-y_{0}+y_{2} & 2y_{1}+y_{3} & -y_{1} & -y_{0}+y_{3}+y_{4} & -y_{2}
\end{array}
\right).$$ \vskip 0.5cm

 If we fix a point $p\in\p^4$ and solve the
system of linear equations $$B(p)\cdot \mathbf{x}^t=0$$ we get the
fiber $\phi^{-1}(p)$. In particular, it is immediate to see that
$\phi^{-1}(0:0:0:0:1)$ is the plane $x_3=x_4=0$ which cuts $X^1$
along a smooth quartic plane curve. An easy computation on
$B(\mathbf{y})$ shows that $\Sing(X^2)=(0:0:0:0:1)$ so that the
are no other 2-dimensional exceptional fibers.

We remark that the general case of this quarto-quartic
transformation is studied in \cite{Todd}, while the existence of a
surface $X^1$ of the above type containing plane quartic curves
is suggested by Room in example iii) of page 80 of \cite{Room}.
Room does not give an explicit example of such an $X^1$
containing plane quartic curves. For this reason we called the
above example(s) {\it Todd--Room contraction}.} \label{surfp4}
\end{example}
\medskip

\begin{example} {\rm (Bordiga conic-bundle). {\it A Fano-Mori
 contraction of fiber type between a smooth projective
4-fold and $\p^3$ with 10 isolated 2-dimensional fibers}. Take
$m=4$ and $n=3$ in the above construction and let us consider a general
matrix of linear forms $A(\mathbf{x})\in M_{4\times
3}(k[\mathbf{x}])$, $k$ an arbitrary field. If it defines a smooth
surface  $X^1\subset\p^4$, then $X^1$ is a degree 6 and genus 3
Bordiga surface, which is arithmetically Cohen-Macaulay. By Bertini's
theorem in characteristic zero and over an algebraically closed
field a general $A(\mathbf{x})$ will define such a surface. A
random choice over a suitable field will also give such examples.

Let $\phi:\p^4\map\p^3$ be the associated map. By the general
discussion in remark \ref{re:secant} $\Sec_3(X^1)=\p^4$ and
through the general point of $\p^4$ there passes a unique
trisecant line to $X^1$, i.e. the Bordiga surface has {\it one
apparent triple point}, see \cite{Sev}. By proposition
\ref{determinantal} $\phi$ is dominant and for general
$\mathbf{y}\in\p^3$ we have $\rk(B(\mathbf{y}))=3$. If
$\rk(B(\mathbf{y}))\leq 2$ for some $\mathbf{y}\in\p^3$, then
necessarily $\rk(B(\mathbf{y}))=2$ and the corresponding plane
$\p^2_{\mathbf{y}}$ cuts $X^1$ along a determinantal cubic plane
curve.  To compute the number of these planes, we can calculate
the degree of the locus of points $\mathbf{y}\in\p^3$ such that
$\rk(B(\mathbf{y}))\leq 2$.  If we assume that $B(\mathbf{y})$ is
a general $3\times 5$ matrix of linear forms on $\p^3$ then the
above locus consists of exactly 10 points by well known formulas,
see for example \cite{Room}, p. 42. This gives the existence of 10
exceptional two dimensional fibers.

We remark that the above description could be also obtained by
looking at  the plane representation of a general Bordiga surface,
eventually passing to $\overline{k}$. One easily verifies that
there are exactly 10 cubics on such a surface, forming a double 10
with the ten chosen points for the plane representation. In fact a
general Bordiga  surface is the blow-up of $\p^2$ at 10 general
points which is embedded by the linear system of plane quartics
through these points.

By proposition \ref{determinantal} a quartic hypersurface through
a Bordiga surface has equation given by the determinant of a
$4\times 4$ matrix of linear forms on $\p^4$ and it is rational
over $k$.

The study of this examples is very classical and attracted the
attention of classical geometers for its remarkable geometry
(congruence of lines, etc, etc.).}\label{Bordiga}
\end{example}
\medskip

\begin{example} {\rm {\it A Fano-Mori birational contraction between a smooth
projective 5-fold and $\p^5$ having one dimensional exceptional
fibers except for an isolated 2 dimensional (or 3 dimensional)
fiber.}

Take $m=n=5$ in the general construction and let
$X^1\subset\p^5$ be a degree 15 and genus 26 smooth 3-fold
defined by a matrix of linear forms $A(\mathbf{x})\in M_{6\times
5}(k[\mathbf{x}])$ for a suitable field $k$.

Since $\Sec_5(X^1)\subsetneqq\p^5$ by remark \ref{re:secant}, we
have a special Cremona transformation $\phi:\p^5\map\p^5$. By
lemma 1.1 $\Sec_5(X^1)$ is a geometrically irreducible
hypersurface (of degree 24). The inverse of $\phi$  is not defined
along  a codimension 2 determinantal scheme $X^2$, which is
generically smooth and defined by the maximal minors of
$B(\mathbf{y})$. The exceptional locus $E_2$ of
$\widetilde{\phi}:Bl_{X_1}(\p^5)=\widetilde{\p^5}\to\p^5$ is then
irreducible of dimension 4 and it is smooth outside
$\widetilde{\phi}^{-1}(\Sing(X^2))$.

From the above geometrical description it follows that, for a
general choice of $A(\mathbf{x})$, $X^2$ is smooth. In fact, for
example by counting parameters we have that a general $X^1$ does
not contain any plane quintic determinantal curve, neither a
quintic determinantal surface spanning a $\p^3$.

Now we  furnish examples of matrices $A(\mathbf{x})$, defined over
suitable fields $k$, such that $X^2$ has only one singular point
$p$ and such that $\widetilde{\phi}^{-1}(p)\simeq\p^2$,
respectively $\widetilde{\phi}^{-1}(p)\simeq\p^3$. The
singularities of both $X_2$ and $E_2$ can be described explicitly
as in the previous example, but we leave the details to the
interested reader.

Let us construct an example such that $X^2$ has one singular point
such that the corresponding fibre of $\widetilde{\phi}$ will be
isomorphic to an "isolated" $\p^2$.
 We can
take $A(\mathbf{x})$ equal to
\medskip

\hskip -1cm
 $\left(
\begin{array}{ccccc}
-2x_{1}-2x_{4} & -2x_{1}+x_{5}-x_{0} & -x_{2}+x_{5}-x_{0} &
x_{3}-2x_{4}+2x_{0} & -2x_{2}+2x_{4} \\
2x_{2}+x_{0} & 2x_{2}-2x_{0} & x_{2}-x_{3}-x_{4} & -x_{1}+x_{2} &
-2x_{3}+2x_{4} \\
x_{4}+x_{5} & -x_{4}+x_{5}+x_{0} & x_{1}+2x_{0} &
-x_{2}-2x_{5}+x_{0} &
x_{1}+x_{5} \\
-x_{3}-x_{0} & -2x_{1}+x_{3} & -2x_{1}+x_{4} & -x_{2}+x_{4} &
-x_{1}-x_{3}-2x_{5} \\
x_{1}+x_{3}-2x_{5} & x_{3}-2x_{5} & -2x_{3}+2x_{5} & x_{1}+2x_{5}
&
x_{2}-x_{0} \\
x_{4} & x_{5} & x_{4} & x_{5} & x_{4}
\end{array}
\right).$
\medskip

To have an example with an isolated exceptional fiber of dimension
3, in fact isomorphic to a $\p^3$, it is sufficient to take
$A(\mathbf{x})$ equal to
\medskip

\hskip -1cm
$\left(
\begin{array}{ccccc}
-2x_{1}-2x_{4} & -2x_{1}+x_{5}-x_{0} & -x_{2}+x_{5}-x_{0} &
x_{3}-2x_{4}+2x_{0} & -2x_{2}+2x_{4} \\
2x_{2}+x_{0} & 2x_{2}-2x_{0} & x_{2}-x_{3}-x_{4} & -x_{1}+x_{2} &
-2x_{3}+2x_{4} \\
x_{4}+x_{5} & -x_{4}+x_{5}+x_{0} & x_{1}+2x_{0} &
-x_{2}-2x_{5}+x_{0} &
x_{1}+x_{5} \\
-x_{3}-x_{0} & -2x_{1}+x_{3} & -2x_{1}+x_{4} & -x_{2}+x_{4} &
-x_{1}-x_{3}-2x_{5} \\
x_{1}+x_{3}-2x_{5} & x_{3}-2x_{5} & -2x_{3}+2x_{5} & x_{1}+2x_{5}
&
x_{2}-x_{0} \\
x_{4} & x_{5} & x_{3} & x_{5} & x_{4}
\end{array}
\right).$
\medskip

It is clear that  the plane $x_3=x_4=x_5=0$, respectively the
3-plane of equation $x_4=x_5=0$, is mapped to $(0:0:0:0:0:1)$ and
it is easy to verify that this is the only singular point of
$X_2$, i.e. that the corresponding linear space is the unique
exceptional fiber of the unexpected dimension}.\label{birp5}
\end{example}
\medskip

We now discuss the last example of this kind, a smooth 3-fold
$X^1\subset\p^5$ of degree 10 and genus 11   having one apparent
quadruple point. Its general hyperplane section is the surface
considered in example \ref{surfp4}. Firstly we discuss the general
case and then we study an exceptional phenomenon for the fibers
\medskip

\begin{example}{\rm(Determinantal 5-dimensional conic bundle) {\it A Fano-Mori
 birational contraction between a
smooth 5-fold and $\p^4$ having an isolated 3-dimensional fiber
in a one dimensional locus of 2-dimensional fibers}. Take $m=5$
and $n=4$ in the general construction. Let us discuss firstly the
general case to construct a conic bundle
$\widetilde{\phi}:\widetilde{\p^5}\to\p^4$ such that the general
fiber is a 4-secant line to $X^1$ and having a smooth curve
$C\subset\p^4$ of 2 dimensional fibers isomorphic to $\p^2$. To
this aim let us take $A(\mathbf{x})$ equal to
\medskip

 $$\left(
\begin{array}{cccc}
x_{1}+x_{4} & x_{1}+x_{5}-x_{0} & -x_{2}+x_{5}-x_{0} & x_{3}+x_{4}-x_{0} \\
-x_{2}+x_{0} & -x_{2}+x_{0} & x_{2}-x_{3}-x_{4} & -x_{1}+x_{2} \\
x_{4}+x_{5} & -x_{4}+x_{5}+x_{0} & x_{1}-x_{0} & -x_{2}+x_{5}+x_{0} \\
-x_{3}-x_{0} & x_{1}+x_{3} & x_{1}+x_{4} & -x_{2}+x_{4} \\
x_{1}+x_{3}+x_{5} & x_{3}+x_{5} & x_{3}-x_{5} & x_{1}-x_{5}
\end{array}
\right).$$
\medskip
The corresponding matrix $B(\mathbf{y})$ is equal to
\medskip

 $$\left(
\begin{array}{cccccc}
y_{0}+y_{4} & -y_{1} & -y_{3}+y_{4} & y_{0}+y_{2} & y_{2}+y_{4} &
y_{1}-y_{3}
\\
y_{0}+y_{3} & -y_{1} & y_{3}+y_{4} & -y_{2} & y_{0}+y_{2}+y_{4} &
-y_{0}+y_{1}+y_{2} \\
y_{2}+y_{3} & -y_{0}+y_{1} & -y_{1}+y_{4} & -y_{1}+y_{3} &
y_{0}-y_{4} &
-y_{0}-y_{2} \\
-y_{1}+y_{4} & y_{1}-y_{2}-y_{3} & y_{0} & y_{0}+y_{3} &
y_{2}-y_{4} & -y_{0}+y_{2}
\end{array}
\right).$$
\medskip

The locus of points $\mathbf{y}$ of $\p^4$ such that
$\rk(B(\mathbf{y}))=3$ is a smooth curve $C$ of degree 20 and
genus 26.

Now we can  force the existence of an exceptional isolated fiber of greater
dimension by taking the following matrix $A(\mathbf{x})$ defining
a smooth 3-fold $X^1\subset\p^5$ as above.
\medskip
$$\left(
\begin{array}{cccc}
x_{1}+x_{2}+x_{4}+x_{5} & x_{1}-x_{5}-x_{0} & -x_{2}+x_{5}-x_{0} &
x_{3}+x_{4}-x_{0} \\
-x_{2}+x_{0} & -x_{2}+x_{3}+x_{0} & x_{2}-x_{3}-x_{4} &
-x_{1}+x_{2}-x_{5}
\\
x_{1}+x_{3}+x_{4}+x_{5} & -x_{4}+x_{5}+x_{0} & x_{1}+x_{3}-x_{0} &
x_{1}-x_{2}+x_{5}+x_{0} \\
-x_{3}-x_{0} & x_{1}+x_{3} & x_{1}+x_{4}-x_{5} & -x_{2}+x_{4} \\
x_{4} & x_{5} & x_{4} & x_{5}
\end{array}
\right).$$
\medskip

In this case the resulting curve $C\subset\p^4$ will have exactly
one singular point at $(0:0:0:0:1)$; the fibre over $(0:0:0:0:1)$
is  isomorphic to $\p^3$. We leave the details to the reader.

We remark that a general hyperplane section of this smooth 3-fold
$X^1\subset\p^5$, would give a whole class of examples of surfaces
$X^1\subset\p^4$ originating examples of Fano-Mori birational
contractions of the type described in example \ref{surfp4}. In
fact all these surfaces will contain a plane quartic curve, the
intersection of the hyperplane with the quartic surface cut on
$X^1$ by $\widetilde{\phi}^{-1}((0:0:0:0:1))\simeq\p^3$.

By proposition \ref{determinantal} a quintic hypersurface through
a variety $X^1\subset\p^5$ as above has equation given by the
determinant of a $5\times 5$ matrix of linear forms on $\p^5$ and
it is rational over $k$.}\label{conicp5}
\end{example}

\section{More preliminary results}

Let us recall that a set of  homogeneous  forms of fixed degree
$d_1\geq 2$ is said to verify {\it condition $K_{d_1}$}, see
\cite{Ve}, if the trivial (or Koszul) syzygies of these forms are
generated by the linear ones. This is a a weakening of the
condition of having the module of first syzygies generated by the
linear ones and it is sufficient for many applications. Particular
forms of this result were known to Room, see \cite{Room}
proposition 2.9.1 for example, as we explained in the previous
section. About {\it condition $K_{d_1}$} we have the following:

\begin{Proposition}{\rm (\cite{Ve}, \cite{HKS})}
 Let $X^1\subset\p^m$ be a smooth variety of dimension $r_1$
scheme theoretically defined by forms $F_0,\ldots,F_n$ of degree
$d_1$ satisfying condition $K_{d_1}$. Then the rational map
$\phi:\p^m\dasharrow\p^n$ defined by the $F_i$'s is an embedding
off $\Sec_{d_1}(X)$. Moreover the closure of a positive
dimensional fiber of $\phi$ is a linear space in $\p^m$
intersecting $X^1$ in a  hypersurface of degree $d_1$. In
particular, if $\Sec_{d_1}(X_1)\subsetneqq\p^n$ and if $m=n$, then
$\phi$ is a special Cremona transformation. \label{prop:Vermeire}
\end{Proposition}
\medskip

In the sequel we will use the above proposition to study the case
of birational maps between a projective space and a cubic
hypersurface in $\p^5$. About the proof of the proposition we
recall that it relies on the fact that, outside the exceptional
divisor of the blow up of $\p^m$ at $X^1$, the graph
$\overline{\Gamma_{\phi}}\subset\p^m\times\p^n$ coincides with the
scheme given by the equations $y_iF_j-y_jF_i$. By conditions
$K_{d_1}$ it can be described by bihomogeneous forms of degree
$(1,1)$.
\medskip

Now let $Y\subset\p^{2n+1}$, $n\geq 2$, be a smooth  hypersurface
defined over a field $k$. Let $\phi:Y\dasharrow \p^{2n}$ be a
birational map of type $(d_1,d_2)$. This means that $\phi$ is
given by a linear subsystem of $|\O_Y(d_1)|$ and $\phi^{-1}$ by a
linear subsystem of $|\O_{\p^{2n}}(d_2)|$. Let $X$ be the base
locus scheme of $\phi$. We will assume that  that $X\subset Y$ is
smooth and connected, i.e. that $\phi$ is {\it special}.

Let $\pi:\widetilde{Y}=Bl_X(Y)\to Y$ stand for the structural map.
As we can consider $\widetilde{Y}$ contained in
$\p^{2n+1}\times\p^{2n}$,  there is a canonical diagram as
follows:
$$\xymatrix{
\widetilde{Y}\ar[d]_{\pi}\ar[dr]^{{\widetilde{\phi}}}       &             \\
  Y  \ar@{-->}[r]_{\phi} & \p^{2n}              }
$$
where $\widetilde{\phi}$ is the projection onto the second factor.

Let $S=\{u \in\p^{2n}\::\: \dim(\widetilde{\phi}^{-1}(u))\geq
1\}$. Since $\p^{2n}$ is smooth, $E=:\widetilde{\phi}^{-1}(S)$ is
a divisor, which is irreducible because
$Pic(\widetilde{Y})=\mathbb{Z}\oplus\mathbb{Z}$, see for example
\cite{ESB}, proposition 1.3. Let $Z$ denote the base locus scheme
of $\widetilde{\phi}^{-1}$, so that $S=Z_{\red}$. Let
$E_1=\pi^{-1}(X)$ and $E_2=\widetilde{\phi}^{-1}(Z)$. Take $H_1\in
|\pi^*(\O_Y(1))|$, $H_2\in |\widetilde{\phi}^*(\O_{\p^{2n}}(1))|$
and let $T_1$ and $T_2$ stand, respectively, for the strict
transform on $\widetilde{Y}$ of a general line contained in $Y$
and of a general line of $\p^{2n}$.

Then we have the following generalization of lemma
\ref{equations}.

\begin{Lemma}\label{equations2}
Let notations be as above.
\begin{enumerate}
\item[{\rm (i)}] $H_2\sim d_1H_1-E_1$, $H_1\sim d_2H_2-E_2$,
$E_2=(E_2)_{\red}$,
 $T_2\cdot H_1=d_2$ and $E_1\cdot T_2=d_1d_2-1$. Moreover, $E_2\sim (d_1d_2-1)H_1-d_2E_1$.
\item[{\rm (ii)}] $\widetilde{\phi}(E_1)=\mbox{\rm Sec}_{d_2}(Z)$
is a hypersurface of degree $d_1d_2-1$; $\pi(E_2)$ is a divisor in
$|\O_Y(d_1d_2-1)|$, which is the locus of the $d_1$-secant lines
to $X$ contained in $Y$. \item[{\rm (iii)}]
$2+r_1=d_2[2n-1+d_1(r_1+1-2n)]$ so that
$\widetilde{\phi}:\widetilde{Y} \to\p^{2n}$ is a Fano-Mori
contraction. \item[{\rm (iv)}] Suppose $\codim(Z)=2$ and let
$V=\p^{2n}\setminus \{z\in Z\,:\,
\dim(\widetilde{\phi}^{-1}(z))\geq 2\}$.  Then
$\widetilde{\phi}^{-1}(V)$ is the blowing up of $V$ along
$Z\setminus \{z\in Z\,:\, \dim(\widetilde{\phi}^{-1}(z))\geq 2\}$,
which is a smooth variety.
 \item[{\rm (v)}] Suppose $k=\overline{k}$ and
let $U=\p^{2n}\setminus \Sing(Z_{\red})$. Then
$\widetilde{\phi}^{-1}(U)$ is the blowing up of $U$ along
$Z_{\red}\cap U$. In $U$, $Z_{\red}=Z$ so that $Z$ is generically
smooth.

\end{enumerate}
\end{Lemma}

Some results of this type were used in \cite{RS} to classify
quadro-quadric special birational maps from $\p^{2m-2}$ and
$\mathbb{G}(1,m)$. Clearly they can be used for the study of
arbitrary birational map from $\p^q$ onto projectively  normal
smooth varieties with Picard group isomorphic to $\mathbb{Z}$,
generated by the hyperplane section, and containing a large family
of lines.

In the sequel we will need the following geometric result, which
relies strongly on the fact that a del Pezzo quintic surface
$X\subset\p^5$ is a variety with one apparent double point. Let us
recall that a smooth irreducible $n$-dimensional variety
$X\subset\p^{2n+1}$ is said to be a {\it variety with one apparent
double point} if through a general point of $\p^{2n+1}$ there
passes a unique secant line to $X$. This means that by projecting
$X$ from a general point of $\p^{2n+1}$ into $\p^{2n}$ the variety
$X$ acquires only one double point as its singularities.

\begin{Lemma} Let $X\subset\p^5$ a smooth del Pezzo surface of
degree 5 and let $Y\subset\p^5$ be a smooth cubic hypersurface
containing $X$. Then $Y$  contains at most a finite number of
planes spanned by the conics contained in $X$.\label{delpezzo}
\end{Lemma}
\begin{proof} The conics contained in $X$ are divided into 5 pencils, at least over
$\overline{k}$. The union of the planes generated by one pencil of
conics is a rational normal scroll, which is necessarily a Segre
3-fold $\p^1\times\p^2\subset\p^5$, see \cite{Ru} and \cite{AR}.
Hence $X$ is contained in five Segre 3-folds defined over
$\overline{k}$. If $Y$ would contain infinitely many planes spanned by conics contained in $X$, then it
necessarily would contain a Segre 3-fold, but this would force $Y$ to
be singular by the Grothendieck-Lefschetz theorem.
\end{proof}

\section{Fano's extremal elementary contractions from
some smooth cubic hypersurface in $\p^5$ onto $\p^4$ or onto a
smooth quadric $Q\subset\p^5$}

In this section we prove the following results,
which were originally proved by Fano in \cite{Fano1} for general cubic hypersurfaces of the mentioned
types.

\begin{Theorem}{\rm (\cite{Fano1})} Let $Y\subset\p^5$ be a smooth
cubic hypersurface defined over a field $k$ and containing a
quintic smooth del Pezzo surface $X$ defined over $k$. Then the
linear systems of quadrics through $X$ restricted to $Y$ defines
an elementary extremal birational contraction
$\widetilde{\phi}:\widetilde{Y}=Bl_X(Y)\to\p^4$. The inverse map
$\widetilde{\phi}^{-1}:\p^4\map Y$ is given by quartic
hypersurfaces passing through a normal surface $Z\subset\p^4$ of
degree 9 and genus 8. For a general choice of $Y$ through $X$, the
surface $Z$ is smooth, and it is a (non-minimal) K3 surface of
degree 9 and genus 8, the blow-up at 5 points of a general degree
14 and genus 8 K3 surface, and $\widetilde{\phi}$  is the inverse
of the blow-up with center $Z$. The fiber over an isolated
singular point of $Z$, if any, is a  plane  contained in $Y$ and
spanned by a conic $C\subset X\subset\p^5$. There exist smooth $Y$
containing a finite number of planes spanned by conics in $X$;
these $Y$ furnish examples of extremal elementary birational
contractions $\widetilde{\phi}:\widetilde{Y}\to\p^4$ having
isolated 2 dimensional exceptional fibers isomorphic to
$\p^2$.\label{th:Fano}
\end{Theorem}
\medskip

\begin{Remark}{\rm It is well known that smooth cubic hypersurfaces $Y$
in $\p^5$ containing a del Pezzo quintic surface $X$ are Pfaffian
and describe an irreducible divisor in the moduli space of smooth
cubic hypersurfaces, see for example \cite{Be}. This divisor
coincides with the locus of smooth cubic hypersurfces containing a
rational normal scroll $X$ of degree 4, see \cite{Fano1}. In fact
this was the main result and the main motivation of Fano's paper.
The theorems, proved by Fano in the second part of \cite{Fano1},
describe every smooth cubic hypersurface through such $X$'s. The
locus of smooth cubic hypersurfaces containing a plane forms a
divisor too. The smooth cubic hypersurfaces giving Fano-Mori
contractions, not inverse of a blow-up, are particular smooth
cubic hypersurfaces in the intersection of the above divisors in
the moduli space.

The general cubic hypersurfaces through such $X$'s  were later
studied by various authors by using Riemann-Roch, Serre's duality,
Kodaira's vanishing (so restricting to characteristic zero and
algebraically closed fields) and similar or more advanced
techniques, see for example \cite{Tregub} and \cite{Tregub2}. Our
approach will be completely elementary and geometric and will
follow Fano's original argument. We use basics of liaison, the
theorem of Danilov and the above generalizations of the
preliminary results of \cite{ESB}.}
\end{Remark}

\medskip

The proof of theorem \ref{th:Fano} will occupy the rest of this
section. We  recall that, by the same techniques, we could also
prove the following theorem.
\medskip

\begin{Theorem}{\rm (\cite{Fano1})} Let $Y\subset\p^5$ be a smooth
cubic hypersurface defined over a field $k$ and containing a
rational normal scroll $X$ of degree 4  defined over $k$. Then the
linear system of quadrics through $X$, restricted to $Y$, defines
an elementary extremal birational contraction
$\widetilde{\phi}:\widetilde{Y}=Bl_X(Y)\to Q\subset\p^5$, where
$Q$ is a smooth quadric hypersurface. The inverse map
$\widetilde{\phi}^{-1}:Q\map Y$ is given by the restriction to $Q$
of the linear system of cubic hypersurfaces passing through a
normal surface $Z\subset Q$ of degree 10 and genus 7. For a
general choice of $Y$ through $X$, the surface $Z$ is smooth and
it is a (non-minimal) K3 surface of degree 10 and genus 7, the
blow-up in a point of a general degree 14 and genus 8 K3 surface,
and $\widetilde{\phi}$  is the inverse of the blow-up with center
$Z$. The fiber over a singular isolated point of $Z$, if any, is a
plane contained in $Y$ and spanned by a conic $C\subset
X\subset\p^5$. There exist smooth $Y$ containing a finite number
of planes spanned by conics in $X$; these $Y$ furnish examples of
extremal elementary birational contractions
$\widetilde{\phi}:\widetilde{Y}\to Q\subset\p^5$ having isolated 2
dimensional exceptional fibers isomorphic to
$\p^2$.\label{th:Fano2}
\end{Theorem}
\medskip

We will not give the proof of  theorem \ref{th:Fano2} here,
but we will state all the needed lemmas and propositions after
the corresponding (and proved) statements needed for the proof of theorem \ref{th:Fano}.
Let us begin with the following Lemma.
\medskip

\begin{Lemma}
Let $X\subset\p^5$ be a smooth quintic del Pezzo surface defined
over a field $k$. Let $H$ be the class of a hyperplane divisor of
$\p^5$. The linear system $|2H-X|=|H^0(I_X(2))|$ on $\p^5$ has $X$
as base locus scheme and it defines a morphism
$\widetilde{\phi}:Bl_X(\p^5)\to \p^4$ such that every fiber is
either isomorphic to a secant line to $X$ or to a plane spanned by
a conic contained in $X$.

Let $X'\subset\p^5$ be a smooth rational normal scroll of degree 4
defined over a field $k$. The linear system $|2H-X'|=|H^0(I_
X'(2))|$ on $\p^5$ has $X'$ as base locus scheme and it defines a
morphism $\widetilde{\phi}':Bl_ {X'}(\p^5)\to \p^5$, whose image
is a smooth quadric hypersurface $Q$, such that every fiber is
either isomorphic to a secant line to $X'$ or to a plane spanned
by a conic contained in $X'$.

\end{Lemma}
\begin{proof} It is well known that the ideal of $X$ is generated
by 5 quadratic forms having the first syzygies generated by the
linear ones. Moreover, by applying proposition \ref{prop:Vermeire}
to the general hyperplane section of $X$, we have that
$\widetilde{\phi}$ is surjective and the general fiber of
$\widetilde{\phi}$ is a secant line to $X$, so that $X$ has one
apparent double point (see also \cite{AR}). The conclusion now
follows taking into account that all the fibers of
$\widetilde{\phi}$ are positive dimensional.
\end{proof}

Let now $Y\subset\p^5$ be a smooth cubic hypersurface defined
over $k$ and  containing $X$. Let $Y'\subset\p^5$ be a smooth cubic hypersurface defined
over $k$ and  containing $X'$,  Let notations be
as in the previous section and by abuse of language let us also
denote by $\widetilde{\phi}:\widetilde{Y}=Bl_X(Y)\to\p^4$  the
restriction of $\widetilde{\phi}$ to $\widetilde{Y}$,
respectively by $\widetilde{\phi}':\widetilde{Y}'=Bl_{X'}(Y)\to
Q\subset\p^5$ the restriction of $\widetilde{\phi}'$ to
$\widetilde{Y}'$.

\begin{Lemma} The morphism
$\widetilde{\phi}:\widetilde{Y}=Bl_X(Y)\to\p^4$ is birational and
its exceptional fibers are isomorphic either to a secant line to
$X$ contained in $Y$, general case, or to a plane contained in $Y$
and spanned by a conic contained in $X$.  The base locus $Z$ of
$\widetilde{\phi}^{-1}$ is a 2-dimensional irreducible scheme with
at most a finite number of singular points, which are exactly the
images of  the planes contained in $Y$, spanned by conics
contained in $X$.

The morphism
$\widetilde{\phi}':\widetilde{Y}=Bl_{X'}(Y)\to Q\subset\p^5$ is birational
and its exceptional fibers are isomorphic either to a
secant line to $X'$, general case, or to a plane spanned by a
conic contained in $X'$.  The base locus $Z'$ of
$\widetilde{\phi}'^{-1}$ is a 2-dimensional irreducible scheme
with at most a finite number of singular points, which are
exactly the images of  the planes contained in $Y'$, spanned by
conics contained in $X'$.

\end{Lemma}

\begin{proof} By the analysis of the fibers of
$\widetilde{\phi}:\widetilde{\p^5}\to\p^4$ given in the previous
lemma, it follows that $\widetilde{\phi}:\widetilde{Y}\to\p^4$ is
a birational morphism: a general secant line to $X$ is not
contained in $Y$ and cuts $Y$ only in one point outside $X$. Hence
the general fiber of $\widetilde{\phi}$ consists of one point. The
other fibers can be only lines (in
$\p^5\times\widetilde{\phi}^{-1}(p)$, the residual intersection of
the plane of a conic in $X$, outside $X$) or eventually planes
(these are the planes of the conics contained in $X$ which are
also planes contained in $Y$). By lemma \ref{delpezzo} there are
only a finite number of these planes contained in $Y$ and the
conclusion follows from lemma \ref{equations2} part iv). In fact
the general fiber of $\widetilde{\phi}_{|E_1}:E_1\to Z$ is a line,
so that that $Z$ is 2-dimensional. \end{proof}

\begin{Remark}{\rm The birationality of the restriction
of $\widetilde{\phi}$ to $\widetilde{Y}$ has the following
algebraic interpretation. The graph of $\phi:\p^5\map\p^4$ is
$\widetilde{\p^5}=Bl_X(\p^5)\subset\p^5\times\p^4$, which is given
by forms of bidegree (1,1), outside the exceptional divisor $E$,
by proposition \ref{prop:Vermeire}. Points
$(\mathbf{x},\mathbf{y})\in(\widetilde{\p^5}\setminus E)
\times\p^4$ are of the following form
$\mathbf{y}=(Q_0(\mathbf{x}):\ldots:Q_4(\mathbf{x}))$, where the
$Q_i$ are the quadratic forms defining $X$. Let $F(\mathbf{x})=0$
be the equation of $Y$. The graph of $\phi:Y\map\p^4$ is
$\widetilde{Y}=Bl_X(Y)$. Outside the exceptional divisor this
graph can be described by the above (1,1) forms together with the
equation of $Y$ "restricted" to $\widetilde{\p^5}\setminus E$.
Since $F=\sum L_iQ_i$, with $L_i(\mathbf{x})$ linear forms, the
restriction of $ F(\mathbf{x})$ to $\widetilde{\p^5}$ outside $E$
is defined by one more (1,1) form: $\sum L_i(\mathbf{x})y_i=0$; so
that the general fiber of $\widetilde{\phi}$ becomes a point. This
means that the cubic hypersurface gives rise to a linear syzygy
among the equations defining $X$.}
\end{Remark}

By the previous lemmas it follows that the 2-dimensional scheme
$Z$ parameterizes the secant lines to $X$ contained in $Y$. Hence
these secant lines form an irreducible 2-dimensional family. The
inverse map $\widetilde{\phi}^{-1}:\p^4\map Y$ will be given by
forms of degree $d_2$ having the base locus supported on $Z$.

\begin{Proposition}{\rm (\cite{Fano1})} Let notations be as above.
Then $\widetilde{\phi}^{-1}:\p^4\map Y$ is given by quartics
vanishing on $Z$, i.e. $d_2=4$, and $\pi(E_2)\in|\O_Y(7)|$.
Moreover, $\widetilde{\phi}$ is the blow-down of $E_2$ outside
$\widetilde{\phi}^{-1}(\Sing(Z))$ and  the quintic del Pezzo
surface is a locus of singular points of multiplicity 4 for
$\pi(E_2)$.

The map $\widetilde{\phi}'^{-1}:Q\map Y'$ is given by cubics
vanishing on $Z'$, i.e. by $|\O_Q(3)-Z'|$, and
$\pi(E_2)\in|\O_{Y'}(5)|$. Moreover, $\widetilde{\phi}'$ is the
blow-down of $E_2$ outside $\widetilde{\phi}'^{-1}(\Sing(Z'))$ and
the rational normal scroll $X'$ is a locus of singular points of
multiplicity 3 for $\pi(E_2)$.

\label{d2}
\end{Proposition}

\begin{proof} Since $r_1=2$ and $n=2$ by lemma \ref{equations2}
part iii) we get $d_2=4$ and that, by lemma \ref{equations2} part
i), $E_2\sim 7H_1-4E_1$, i.e. {\it le $\infty^2$ rette del sistema
$\Gamma$ ricoprono una $M^{7\cdot 3}$ per la quale $X$ \' e
superficie quadrupla}, see \cite{Fano1}, pg. 79.

We cannot resists in furnishing Fano's original geometric argument
on page 79 of \cite{Fano1} proving that $d_2=4$, also to show the
richness of his geometrical arguments and  how it works perfectly.
Let $\Gamma:=\pi(E_2)$ be the 2-dimensional family of secant lines
to $X$ contained in $Y$. Fano shows that through a general point
$p\in X$ there passes 4 lines of $\Gamma$, i.e. if $\Gamma\sim
\alpha H_1-\beta E_1$, then $\beta=4$ and apparently he knows that
$\beta=d_2$. Indeed, $E_2\cdot \pi^{-1}(p)=(\alpha H_1-\beta
E_1)\cdot \pi^{-1}(p)=-\beta E_1\cdot \pi^{-1}(p)=\beta$. The
$\beta$ lines of $\Gamma$ passing through $p$ are contained in
$T_p(Y)$. Let $C=T_p(Y)\cap X$. Obviously $C$ has degree 5 and it
is a tangent hyperplane section of $X$, so that it has a double
point in $p$. If we project $C$ from $p$ to a general $\p^3$
inside $T_p(Y)$ we get a 2-dimensional cone over a twisted cubic.
The lines of $\Gamma$ passing through $p$ are the intersections of
$Y$ with this cone, out of $C$. As the intersection scheme has
degree 9 we get: $\beta=9-5=4$. From this he deduces $d_2=4$

Perhaps Fano was not aware that $d_2=4$ implies
$\Gamma\in|\O_Y(7)|$, so he shows this last fact geometrically in
the following way, giving another different proof of the above
assertion, since by lemma \ref{equations2} part i) we have
$7=2d_2-1$. Let us recall that the 3-dimensional $\Gamma$ is the
family of the secant lines to $X$ contained in $Y$. If
$\Gamma\in|\O_Y(d)|$, then $\Gamma$ intersects a general line
$l\subset Y$ in $d$ points. Then $d$ is also the number of lines
of $\Gamma$ intersecting $l$ since through a general point of $Y$
there passes a unique secant line to $X$. By projecting $X$ from
$l$ onto a disjoint $\p^3$ we get a
 degree 5 singular surface $\widetilde{X}:= \pi_l(X)$, as $l$ does not
intersect $X$. One immediately sees that the double curve $D$ of
$\widetilde{X}$ has degree 5, since a general plane section of
$\widetilde{X}$ is a plane quintic of geometric genus 1, the same
sectional genus of $X$, so that it has 5 double points. Now,
 as $X$ has one apparent double point, we can define $\Sigma_l$, the scroll of secant lines to $X$
 passing through (a general point of) $l$.
 A general $\p^4$ through $l$ cuts $\Sigma_l$ along a curve
 which is the union of $l$ and the 5 secant lines to $X$
 which are contracted by $\pi_l$, hence $\deg(\Sigma_l)=6$.
 The secant lines of $\Gamma$ intersecting $l$ are the
 intersection of $Y$ and $\Sigma_l$ outside $l$ and outside
the curve $\pi_l^{-1}(D)$.
 Note that $\pi_l^{-1}(D)$ has degree 10
 because $\pi_l$ is a degree 2 map onto $D$. Therefore $d=6\cdot 3-1-10=7$, see \cite{Fano1}, pg. 79.
\end{proof}
\medskip

The following lemma will conclude the proof of  theorem
\ref{th:Fano}.
\medskip

\begin{Lemma} The base locus scheme $Z\subset\p^4$ of
$\widetilde{\phi}^{-1}$ is a degree 9 and sectional genus 8 normal
surface containing 5 skew lines. The ideal sheaf of $Z$ has the
following minimal resolution: $$0\to\O(-6)\to\O(-5)^{\oplus
6}\to\O(-4)^{\oplus 6}\to \I_Z\to 0.$$ If $Z$ is smooth, then the
5 skew lines are  (-1)-curve and $Z$ is  the blow-up at 5 points
of a genus 8 and degree 14 $K3$ surface
$\widetilde{Z}\subset\p^8$; moreover, the surface $Z$ is  a linear
projection of $\widetilde{Z}$ by a suitable 3-plane cutting
$\widetilde{Z}$ in 5 points. If $Z$ is singular, its singular
points lie on these skew lines, they are of multiplicity 3 and
their tangent cones are cones over a twisted cubic curve; there
are at most 2 singular points on each of the 5 lines. In
particular there are at most 10 singular points  on $Z$.

The base locus scheme $Z'\subset Q\subset\p^5$ of
$\widetilde{\phi}'^{-1}$ is a degree 10 and sectional genus 7
normal surface containing a conic $\gamma$. If $Z'$ is smooth,
then the conic $\gamma$ is a  (-1)-curve and $Z'$ is  the blow-up
at a point of a genus 8 and degree 14 $K3$ surface
$\widetilde{Z}\subset\p^8$; moreover, the surface $Z$ is the
projection of $\widetilde{Z}$ from a general tangent plane. If
$Z'$ is singular, its singular points lie on $\gamma$; they are of
multiplicity 3 and their tangent cones are cones over a twisted
cubic curve.
\end{Lemma}

\begin{proof}
By the previous lemma we know that $Z$ has at most a finite number
of singular points and we calculate its degree. Let us recall that
$\widetilde{\phi}$ is given by the linear system $|2H_1-E_1|$ and
let us take two general divisors $D_1, D_2\in |2H_1-E_1|$. Then
$D_1\cap D_2:=F$ is a smooth surface of degree 7 and sectional
genus 3. By restricting the map $\widetilde{\phi}$ to $F$, we
obtain a morphism $\widetilde{\phi}:F\to \p^2$. Its inverse is
given by the linear system of plane quartic curves passing through
$\deg(Z)$ points by the explicit description of
$\widetilde{\phi}^{-1}$. These points are smooth points for $Z$,
otherwise the linear system of plane quartics would have some
fixed singular point and $F$ could not have sectional genus 3  (by
the way this gives a direct proof that $Z$ is generically reduced
and hence generically smooth). Since $F$ has degree 7, we have
that $\deg(Z)= 4\cdot4-7=9$ (see the footnote 25 on page 79 of
\cite{Fano1}). Moreover lemma \ref{equations2} yields that
$\widetilde{\phi}(E_1)=\Sec_4(Z)$ is a hypersurface of degree 7.

Let us now show that $Z$ is locally Cohen-Macaulay.  To get the
resolution of the ideal sheaf of $Z$ we will use some basic facts
of liaison.

Let us take two general elements in $|4H_2-E_2|$, they give rise
to two general quartic hypersurfaces $V_1$ and $V_2$ in $\p^4$,
containing $Z$, such that $V_1\cap V_2=Z\cup Z"$ as schemes, with
$Z"$ a degree 7 subscheme of $\p^4$. Let us recall that, by lemma
\ref{equations2}, we know that $|4H_2-E_2|$ is the linear system
on $\widetilde Y$ corresponding to the hyperplane sections of $Y$.
Hence, via $\widetilde{\phi}^{-1}$, $Z"$ will correspond to the
intersection of $Y$ with a general $\p^3$; i.e. $Z"$ is the image
by $\widetilde{\phi}$ of the blow-up at 5 $(=\deg(X))$ points of a
smooth cubic surface $G$ obtained by intersecting $Y$ with the
$\p^3$ corresponding to $V_1\cap V_2$. The surface $G$ is smooth
by the general choice of $V_1$ and $V_2$. The restriction of
$\widetilde{\phi}$ to $G$, $\widetilde{\phi}:G\to Z"$, is given by
the linear system of quadric surfaces in $\p^3$ passing through
the five points of intersection of the above $\p^3$ with $X$. By
this description we get that $Z"$ is reduced, irreducible and
normal. We can also suppose that $Z"$ is smooth at $\Sing(Z)$
because a general $\p^3$ will cut any plane
$\widetilde{\phi}^{-1}(z)$, $z\in \Sing(Z)$, at one point only. In
any case $Z"$ is locally Cohen-Macaulay and normal, so that $Z$ is
locally Cohen-Macaulay too by liaison and since it is non-singular
in codimension 1 it is also normal. In conclusion $Z$ is a
reduced, normal, irreducible surface of degree 9.

If $\Sing(Z)$ is not empty a general quartic hypersurface $V$ in
$\p^4$ passing through $z\in \Sing(Z)$ will have a double point at
$z$. In fact the multiplicity of such a $V$ at $z\in \Sing(Z)$
will be the multiplicity at $z$ of the quartic curve obtained by
cutting $V$ with a general plane of $\p^4$ passing through $z$.
Let us consider the restriction of $\widetilde{\phi}^{-1}$ to a
generic plane of $\p^4$ passing through $z$. The image of this
plane, by the restriction of $\widetilde{\phi}^{-1}$, is a surface
in $\p^5$ which is linked $(4,4)$ to the union of $X$ and the
corresponding plane over $z$, so that it has degree 6.  The map is
given by the linear systems of plane quartic passing through $z$
with multiplicity $h$ and other $t$ points with multiplicity 1
such that $10=h^2 +t$. As the projective dimension of that system
of plane quartics is 5, it is easy to see that $h=2$ and $t=6$ so
that $z$ is a point of multiplicity $\deg(Z)-t$ = 3 for $Z$.

 Now, let $V_1, V_2$ be two generic quartic hypersurfaces in $\p^4$ passing through $Z$,
 such that $V_1\cap V_2=Z\cup Z"$. By the above discussion we have
 that $Z"$ is smooth at any  $z\in \Sing(Z)$.
 The tangent cone of $Z$ at any singular point $z$ is a cone over a twisted cubic,
in fact $V_1$ and $V_2$ have multiplicity 2 at $z$ and their
tangent cones at $z$ are cones over quadrics of $\p^3$, so that
the tangent cone  of $Z$ at $z$ is a degree 3 cone over a degree 3
curve which is contained in the intersection of two smooth
quadrics and residual to a line, i.e. a twisted cubic. By a
generalization of \cite{PS}, proposition 4.1 (see also remark 1.12
of \cite{ARa}), we deduce that, for general $V_1, V_2$, $Z"$ is
smooth.

In conclusion the surface $Z$ is linked (4,4) to a smooth surface
$Z"$ having degree 7 and sectional genus 4; for this last fact
recall that $Z"$ is birational to the above mentioned $G$ via the
restriction of $\widetilde{\phi}$. By looking at the resolution of
the ideal sheaf of such a surface in \cite{DES} and by using
liaison we get that the ideal sheaf of $Z$ has minimal resolution
$$0\to\O(-6)\to\O(-5)^{\oplus
6}\to\O(-4)^{\oplus 6}\to \I_Z\to 0.$$

 Let us assume that $Z$ is smooth. By liaison we also get that the sectional genus of $Z$ is 8,
that $q(Z)=h^1(\O_Z)=0$ and that $h^0(\O_Z(K_Z))=1$, since
$h^0(\I_{Z"}(3))=1$. Let us recall that $Z$ contains five skew
lines $l_1,...,l_5$, corresponding to the five pencils of conics
contained in $X$ (see lemma  \ref{delpezzo}). Let $H$ be a
hyperplane section of $Z$ so that $14=H^2+H\cdot K_Z=9+H\cdot K_Z$
implies $H\cdot K_Z=5$. One gets $(K_Z-\sum l_i)\cdot H=0$ so that
$K_Z\sim \sum l_i$ and the 5 lines are the only (-1)-curves of
$Z$.  By blowing down these 5 exceptional (-1)-curves, we obtain a
degree 14 and genus 8 smooth surface $\widetilde{Z}$  such that
$K_{\widetilde{Z}}$ is trivial and $q(\widetilde{Z})=0$, i.e. a
genus 8 smooth $K3$ surface.

 Let us assume that $Z$ is singular. We know that the singular points of $Z$ lie on the five lines quoted above.
  By lemma \ref{delpezzo} we also know that $X$ is a divisor of type (1,2)
  in five Segre 3-folds of $\p^5$ and that the singular points of $Z$
  lying on each line are in one to one correspondence with the planes
  of type (1,0) contained in the intersection of each Segre 3-fold
  with the cubic hypersurface Y. This hypersurface cuts a divisor of type (3,3)
  on each Segre 3-fold. As $X\subset Y$ we have a residual divisor of type (2,1)
  on each Segre 3-fold, and such divisors contain at most two planes of type (1,0).
  So that there are at most 10 singular points on $Z$.
\end{proof}

 The previous lemma concludes the proof of theorem \ref{th:Fano}.
 We can also analyze the singularities of $E_2$ along a plane $\Pi$,
corresponding to $z\in \Sing(Z)$, respectively $z'\in \Sing(Z')$,
\cite{Fano1}, pg. 79.

\begin{Proposition} Let notations be as above. Then the general
point of $\widetilde{\phi}^{-1}(z)$, $z\in \Sing(Z)$, is a point
of multiplicity 2 for $E_2$.

 The general point of $\widetilde{\phi}'^{-1}(z')$, $z'\in \Sing(Z')$,
 is a point of multiplicity 2 for $E_2$.
\end{Proposition}

\begin{proof} Let notations be as in the proof of proposition
\ref{d2}. Let $\Pi$ be any plane, contained in $Y$, and spanned by
a conic $C$ contained in $X$. Let us take a general line $l$ in
$Y$ and let $q$ be the intersection of $\Pi$ with $l$. The degree
6 scroll $\Sigma_l$ described in proposition \ref{d2} degenerates
into the union of $\Pi$ and a degree 5 scroll $\Sigma'_l$ with the
same directrix $l$. Since we can suppose that through every point
of $l\setminus q$ there passes a unique secant line to $X$, the
scroll $\Sigma'_l$ is smooth outside $\Pi\cap \Sigma'_l$ and there
passes only one line of this scroll through $q$. The double curve
$D$ degenerates into the union of a quartic curve and a line,
which is the projection of the conic $C=\Pi\cap X$ from $l$ (see
the proof of proposition \ref{d2}). The intersection with
$\Sigma'_l$ of a general hyperplane of $\p^5$ passing through $l$
consists of four lines outside $\Pi$ and a secant line $l'$ to $C$
through $q$. Then the number of lines of $\Gamma=\pi(E_2)$
intersecting $l$ outside $q$ is equal to $3\cdot 5-8-1-1=5$ where
8 is the degree of the curve projecting itself onto $D$ from $l$
and the last two  -1 are given by the existence of $l$ and $l'$.
This finally yields  $\mult_q(E_2)=2$ because
$\Gamma\in|\O_Y(7)|$.
\end{proof}

\medskip

The fact that there exist smooth cubic hypersurfaces $Y$
containing a smooth del Pezzo surface and some planes of the
conics on the surface is part of the construction of a del Pezzo
surface, or of a rational normal scroll of degree 4.

\begin{Proposition} There exist smooth cubic hypersurfaces $Y\subset\p^5$
containing a smooth quintic del Pezzo surface $X$ and some of the
planes spanned by a conic contained in $X$, respectively
containing a rational normal scroll of degree 4 $X'$ and some of
the planes spanned by conics contained in $X'$.
\end{Proposition}
\begin{proof}
Take the plane $\Pi:\;x_0=x_1=x_2=0$ in $\p^5$ and consider the
Segre variety $M=\p^1\times\p^2$ given by the vanishing of the
$2\times 2$ minors $Q_i$, $i=0,1, 2$ of the matrix of linear forms
$$\left(
 \begin{array}{ccc}
 x_{0}&x_1& x_2\\
 x_{3}&x_{4}&x_{5}
\end{array}
\right).
 $$
Take  a general smooth quadric  hypersurface $Q_3$ containing the
plane $\Pi$. Then $Q_3\cap M=\Pi\cup X$ with $X$ a smooth del
Pezzo surface of degree 5, realized as a divisor of type $(1,2)$
in $M$. Take $Y=l_0Q_0+\ldots+l_3Q_3$ with $l_i=0$ general
hyperplanes in $\p^5$. Then by Bertini theorem, $Y$ is non
singular outside $\Pi\cup X$ and a direct computation shows that
it is also non-singular along $\Pi\cup X$. The ideal of $X$ can be
generated by taking a fifth quadric $Q_4$ through $X$, not
containing $M$ and passing through another general plane of $M$.
From this it is clear that the map $\phi:\p^5\map\p^4$ contracts
the variety $M$ onto a line corresponding, for example, to the
pencil $<Q_3,Q_4>$. The fiber of $\phi_{|M}$ is the unique $\p^2$
of the Segre residual to $X$ obtained by intersecting the
corresponding quadric with $M$. \end{proof}
\medskip

One easily constructs examples of smooth cubic hypersurfaces $Y$
through $X$,
 defined over "sufficiently" small fields and
 containing only one plane spanned by  the conics on $X$, or two planes,
 belonging either to the same pencil or
 to different pencils. We leave the computational details to the
 interested reader.
 As we saw, 10 is the maximal number of planes of the above type
 contained in a smooth cubic hypersurface,
at most  two planes for each pencil, but we do not know if this
maximal case really occurs.

\medskip

We end the section by remarking some amusing consequences in terms
of bilinear forms defining the graph of $\widetilde{\phi}$ inside
$\p^5\times\p^4$, or better $\widetilde{Y}\times \p^4$. We change a little the notations introduced in section 2.

 From the resolution of $\I_Z$ we know that there exists a $6\times 6$ matrix $B(\mathbf{y})$
 of linear forms given by the linear syzygies of  six degree 4
 generators of $I_Z$. By the above description and from the fact that
$\widetilde{\phi}^{-1}$ is given by the linear system of these
quartics, we get some restrictions which $B(\mathbf{y})$ has to
satisfy and we get equations for $Y$. Let us remark that
$\widetilde{Y}=Bl_X(Y)=Bl_Z(\p^4)$ is clearly the graph of
$\widetilde{\phi}$ inside $\p^5\times\p^4$. Outside $E_2$ this
graph is given by the equations $B(\mathbf{y})\cdot
{\mathbf{x}}^t=\mathbf{0}$ which are forms of bidegree $(1,1)$,
where $\mathbf{x}=(x_0,\ldots,x_5)$. Take a (general) point
$\mathbf{p}\in\p^5$. To get the coordinates of
$\widetilde{\phi}(p)$ one has to resolve the system of linear
equations in $\mathbf{y}=(y_0,\ldots,y_4)$ given by
$\mathbf{0}=B(\mathbf{y})\cdot
{\mathbf{p}}^t={A(\mathbf{p})}^t\cdot {\mathbf{y}}^t$. The matrix
$A(\mathbf{x})$ can be obtained, as always, by derivation with
respect to the variables $\mathbf{y}$ of the matrix of homogeneous
$(1,1)$ forms $B(\mathbf{y})\cdot {\mathbf{x}}^t$. We have that
there will be no solutions for a general $\mathbf{p}$ because
$\widetilde{\phi}^{-1}(\p^4)$ is the cubic hypersurface $Y$. This
means that the $6\times 5$ matrix $A(\mathbf{x})$ has rank 5 at
the general point of $\p^5$. On the contrary, at the general point
of $Y$ there will be a unique solution,  translating into the fact
that the six $5\times 5$ minor of $A(\mathbf{x})$ contain a common
irreducible factor of degree 3 which is an equation of $Y$. The
six residual quadratic forms define a (proper) closed subscheme of
$Y$, which in fact has to be a subscheme of $\pi(E_2)=\Gamma$. By
a computer algebra system one easily verifies that this scheme is
the original (or new, depending on the point of view)  quintic del
Pezzo surface.

\section{Semple and Tyrrell extremal elementary contraction  from
a smooth 6-dimensional variety to $\p^6$}

Let $X^1$ be the embedding in $\p^6$ of the blow-up $\tilde \p^2$
 of $\p^2$ at eight points $p_i$, $i=1,\ldots 8,$ in general position,
 given by the linear system of quartics passing through the $p_i'$s,
 (actually it is sufficient that they do not belong to a smooth
 conic and that no four of them are collinear).
 The surface $X^1$ is arithmetically Cohen-Macaulay and
 it is cut out by $7$ quadrics since it is linearly normal, regular
and it has degree $8=2\codim(X^1)$, see \cite{AR}. During the rest
of this section we will prove the following result.
\medskip

\begin{Theorem}{\rm (Semple and Tyrrell birational Fano-Mori
contraction, \cite{ST2})} Let notations be as above. The map
$\phi:\p^6\map\p^6$ given by $|H^0(\I_{X^1}(2))|$ is a special
Cremona transformation of type $(2,4)$. The inverse transformation
is not defined along a subvariety $X^2\subset\p^6$ of dimension 4
and degree 8, which is linked $(4,4)$ to a general projection of
the Veronese embedding $\nu_2(\p^4)$ from the $\p^7$ spanned by 8
general points on it. The variety $X^2$ has a singular locus
consisting of 28 isolated triple points $p_{i,j}$, $i<j$,
$i=1,\ldots, 8$, (corresponding, via $\phi$, to the planes of the
28 conics contained in $X^1$
 which are given by the $\binom{8}{2}$ lines
 joining any couple of the 8 points in $\p^2$)
 and of a quadric surface of double points.
The fibre over any double point of $X^2$ is a quadric surface. The
triple points on $X^2$ have  a cone over a Segre 3-fold
$\p^1\times\p^2$ as tangent cone. The  map $\phi$ gives a
birational extremal elementary contraction
$\widetilde{\phi}:\widetilde{\p^6}=Bl_{X^1}(\p^6)\to\p^6$. The
fibers of $\widetilde{\phi}_{|E_2}\to X^2$ are either proper
transforms of secant lines to $X^1$ or proper transforms of the 28
planes of the conics contained in $X^1$, mapping to the 28
isolated triple points $p_{i,j}$ of $X^2$ or quadric surfaces
mapping to the points of a smooth quadric $\widetilde{Q}$
contained in $X^2$ (these being 2-dimensional non-isolated
fibers). On $\p^6\setminus(\widetilde{Q}\cup(\cup_{i,j} p_{i,j}))$
the morphism $\widetilde{\phi}$ is the blow-up of the smooth
variety $X^2\setminus(\widetilde{Q}\cup(\cup_{i,j}
p_{i,j}))$.\label{SempleTyrrell}
\end{Theorem}
\medskip

Let us briefly show geometrically, as indicated in \cite{ST2}, why
the quadrics through $X^1$ define a Cremona transformation
$\phi:\p^6\map\p^6$ (see \cite{ST2}, \cite{HKS} and also
\cite{AR}, section 6). At the same time we furnish the complete
description of the fibers of $\widetilde{\phi}$.
\medskip

\begin{Lemma} {\rm (\cite{ST2}, \cite{HKS})}
The surface $X^1\subset\p^6$ is the complete intersection of two
divisors of type $(1,2)$ on a cone $W\subset\p^6$ over a Segre
3-fold $\p^1\times\p^2\subset\p^5\subset\p^6$. The vertex of $W$
is $P:=f(p_9)$ where $p_9$ is the ninth  associated point to
$p_1,\ldots,p_8$ and $f:\widetilde{\p^2}\hookrightarrow\p^6$ is
the embedding the surface $\widetilde{\p^2}$.
\end{Lemma}
\begin{proof}

As $f:\widetilde{\p^2}\hookrightarrow\p^6$ is given by the linear
system $|4H-p_1-\ldots -p_8|$, where $H$ is the class of a line in
$\p^2$, each member of the pencil of cubics
$\{D_{\lambda}\}_{\lambda\in\p^1}$ through $p_1,\ldots,p_9$ is
embedded by $f$ as a quartic normal curve, and the general one
will be elliptic. The linear spans $\p^3_{\lambda}=<D_{\lambda}>$,
$ \lambda\in\p^1$, generate a cone $W\subset\p^6$ of vertex $P$
over a Segre 3-fold $\p^1\times\p^2\subset\p^5$, which is the only
smooth (rational) 3-dimensional scroll in $\p^5$. The cone $W$
contains $X^1$ and this surface is the complete intersection of
two divisors of type $(1,2)$ on $W$ (see \cite{ST2}, pg. 214 or
\cite{HKS}, pg. 434), where a divisor $D$ of type $(1,2)$ on $W$
is the intersection of $W$ with a hyperquadric  $Q\subset \p^6$
containing a $\p^3=R$ of the ruling, i.e. $Q\cap W=D\cup R$. Also
from this description one can deduce that the ideal of $X^1$ is
generated by 7 quadratic forms, three coming from the cone over
the Segre 3-fold and the remaining four coming from the fact that
each divisor of type $(1,2)$ on $W$ is given, on $W$, by the
restriction of two quadratic forms,  see \cite{ST2}, pg. 214 or
\cite{HKS}, pg. 434.
\end{proof}
\medskip

\begin{Lemma} {\rm (\cite{ST2})}
For any point $p\in\Sec(X^1)\setminus W$ let $L_p$ be the locus of
secant lines to $X^1$ passing through $p$. Then $L_p$ consists either
of a unique secant (or tangent) line $\p^1_p$ or
of the whole plane $\Pi_{i,j}$ if $p$ belongs to one
of the $28$ planes $\Pi_{i,j}=<C_{i,j}>$ spanned by the 28 conics on $X^1$. \label{locus}
\end{Lemma}

\begin{proof}

Let us take a generic point $p\in\p^6\setminus W=\Sec(W)\setminus W$.
 The locus of secant lines to $W$ passing through $p$ is a linear $\p^4_p$
 cutting $W$ along a rank 4 quadric $Q_p$ which is singular at $P$, the vertex of $W$,
by an obvious adaptation of  example 1.
 Hence $L_p$ is given by the secant lines to the
intersection of $Q_p$ with two divisors of type $(1,2)$ on $Q_p$,
$\widehat{Q_1}$, $\widehat{Q_2}$, i.e. by secant lines through $p$
to a degree 4 curve $C_p\subset \p^4_p$ (see \cite{ST2} pg. 206).
Let us remark that $Q_p\cap\widehat{Q_1}\cap\widehat{Q_2}=C_p\cup
P$. The degree 4 curve $C_p$ can be:
\begin{itemize}
\item[(i)] a smooth quartic rational normal curve $F$, the
embedding by $f$ of the blow up of a line of $\p^2$ not passing
through any of the points $p_i$; \item[(ii)] the union of a {\it
twisted cubic} $T_j$ and a line $E_j$ intersecting at one point,
the embedding by $f$ of the blow up of a line passing through a
point $p_j$; \item[(iii)] the union of a {\it conic} $C_{i,j}$ and
two lines cutting the conic at two distinct points, the embedding
by $f$ of the blow up of a line $<p_i,p_j>$. \end{itemize}
Moreover we have that the locus of secant lines to $C_p$ passing
through a point $p\in\Sec(X^1)\setminus W$ consists either of a
unique secant or tangent line to one of the above reduced curves
or of the whole plane $\Pi_{i,j}$, if $p$ belongs to one of the
$28$ planes $\Pi_{i,j}=<C_{i,j}>$. Recall that a quartic curve in
$\p^4$ has a 4-secant plane only if it is the union of four lines;
 see also the next proposition where it is shown that systems of
quadrics defining these curves satisfy condition $K_2$ yielding
a(nother) proof of this geometrical fact by proposition
\ref{prop:Vermeire}.
\end{proof}

We study the following more general situation.
\begin{Lemma}
Let $C\subset\p^4$ be a reduced non-degenerated quartic curve as in one of the cases
 $(i),(ii),(iii)$ of the previous lemma.
 Suppose also that $C$ is the complete intersection
 of two divisors of type $(1,2)$ of a rank 4 hyperquadric $Q\subset\p^4$,
 whose vertex will be indicated by $O$,  $O\not\in C$.
 Then the rational map defined by the linear system of quadrics
 passing through $C$ defines a birational map $\mu:\p^4\map\widetilde{Q}\subset\p^5$, where
$\widetilde{Q}$ is an irreducible quadric. The subsystem of quadrics
passing also through $O$ gives a birational map
$\mu':\p^4\map\p^4$ not defined along $C\cup O$ and which is the
composition of $\mu$ with the projection of $\widetilde{Q}$ from
$\mu(O)$. The map $\mu'$ is an isomorphism outside $\Sec(C)\cup
\mu^{-1}[T_{\mu(O)}(\widetilde{Q})\cap\widetilde{Q}]=\Sec(C)\cup
Q$.\label{quartic}
\end{Lemma}
\begin{proof} The curve
$C$ is a hyperplane section of one of the two rational scrolls of
degree 4 in $\p^5$, $S(2,2)$ or $S(1,3)$, which are scheme
theoretically defined by 6 quadrics satisfying condition $K_2$ so
that $C$ is scheme theoretically defined by 6 quadrics satisfying
condition $K_2$. The rational map defined by these quadrics is a
birational map, $\mu:\p^4\map\widetilde{Q}\subset\p^5$, where
$\widetilde{Q}$ is an irreducible hyperquadric and $\mu$ is an
isomorphism outside $\Sec(C)$, (see for example \cite{ST2} or
apply proposition \ref{prop:Vermeire} and a simple calculation
with Chern classes to deduce that $\widetilde{Q}$ is a quadric
hypersurface). The subsystem of quadrics passing also through $O$
gives a birational map $\mu':\p^4\map\p^4$ not defined along
$C\cup O$ and which is the composition of $\mu$ with the
projection of $\widetilde{Q}$ from $\mu(O)$. The map $\mu'$ is an
isomorphism outside $\Sec(C)\cup
\mu^{-1}(T_{\mu(O)}(\widetilde{Q})\cap\widetilde{Q})=\Sec(C)\cup
Q$. In fact $\mu^{-1}(T_{\mu(O)}(\widetilde{Q})\cap\widetilde{Q})$
is a quadric hypersurface in $\p^4$ passing through $O$ and $C$,
i.e. an element of the linear system defining $\mu'$. Since the
exceptional locus of $\mu'$ consists only of this quadric
hypersurface and $\Sec(C)$ and, by the definition of $C$, $Q$ is
contracted by $\mu'$, the conclusion easily follows.
\end{proof}

Now we can describe the fibers of $\widetilde{\phi}$ and
we can give a geometrical idea of the birationality of $\phi:\p^6\map\p^6$.

\begin{Lemma} Let notations be as above. The map $\phi:\p^6\map\p^6$
is birational and it is an isomorphism outside $\Sec(X^1)$. Let
$\widetilde{\phi}:\widetilde{\p^6}=Bl_{X^1}\p^6\to\p^6$ its
resolution. Moreover, for any point $p\in\Sec(X^1)\setminus W$ we
have that  $\widetilde{\phi}^{-1}(\widetilde{\phi}(p))$ is equal
to the locus of secant lines to $X^1$ passing through $p$, i.e. is
isomorphic to a secant or tangent line to $X^1$ or to one of the
28 planes $\Pi_{i,j}$ if $p\in\Pi_{i,j}\setminus X^1$. If
$\widetilde{W}$ is the strict transform of $W$ on
$Bl_{X^1}(\p^6)$, then for every $p\in\widetilde{W}$ we have that
$\widetilde{\phi}^{-1}(\widetilde{\phi}(p))$ is isomorphic to  a
quadric surface.
\end{Lemma}

\begin{proof}
We can suppose that,  among the seven quadrics defining $X^1$, $Q_0,
Q_1, Q_2$ define the cone $W$. Let $\psi:\p^6\setminus W\to\p^2$
be the rational map defined by the above net of quadrics.
The closure of every fiber $\psi^{-1}(\psi(p))$, $p\in \p^6\setminus W$ is a linear
space equal to $\p^4_p$, see example 1.
Hence the $\p^4_p$'s vary in  a two dimensional family
such that through every point of $\p^6\setminus W$
there passes only one member of the family.
Then the closure of any fiber $\phi^{-1}(\phi(p))$,
$p\in \p^6\setminus W$, is clearly contained in $\p^4_p$,
so that it is equal to the closure in $\p^4_p$ of
the restriction of $\phi$ to $\p^4_p$. For any $p\in \p^6\setminus W$, the map
$\phi_p:=\phi_{|\p^4_p}$ is not defined along the quartic curve
$C_p$ and in $P$ and it is given by the linear system of quadrics scheme
theoretically defining $C_p\cup P$, see lemmas \ref{locus} and
\ref{quartic}. By lemma \ref{quartic}  the map $\phi_p:\p^4_p\map\p^4\subset\p^6$ is birational so that $\phi$ is
birational.

Moreover,  we can now describe all the positive dimensional fibers
$\phi^{-1}(\phi(p))$ for every $p\in \Sec(X^1)\setminus W$ by
describing the fibers of the corresponding $\phi_p$.

Putting all the above geometric descriptions together we finally
obtain that the closure of the fiber $\phi^{-1}(\phi(p))$, $p\in
\p^6\setminus X^1$ is one of the following types: it reduces to
$p$, if $p\not\in \Sec(X^1)$; it is equal to a tangent line or a
secant line to $C_p$, if $p\in \Sec(X^1)\setminus W$ and $C_p$ has
at most two components; if $C_p$ contains one of the 28 conics
$C_{i,j}$ of $X^1$, then it is a secant or tangent line to $C_p$
if $p\in \Sec(X^1)\setminus ( \Pi_{i,j}\cup W)$ and it is the
whole plane $\Pi_{i,j}$ if $p\in\Pi_{i,j}\setminus C_{i,j}$; it is
a hyperquadric if $p\in W\setminus X^1$, see \cite{ST2}. In fact
the cone $W$ is contracted by $\phi$ onto a smooth quadric surface
$\widetilde{Q} \subset\p(<Q_3,Q_4,Q_5,Q_6>)$. The closure of a
fiber of $\phi_{|W}:W\map \widetilde{Q}$ is easily seen to be a
quadric surface, see also \cite{ST2}.
\end{proof}

Now we are able to conclude the proof of theorem 5.1
by describing the birational map $\widetilde{\phi}:Bl_{X_1}(\p^6)\to\p^6$.

\begin{proof} (of theorem 5.1) By the previous lemma we have that
the strict transform of $\Sec(X^1)$ is 5-dimensional and it is
contracted onto a 4 dimensional variety
$Z=\phi(Sec(X^1))=(X^2)_{\red}$. The planes $\Pi_{i,j}$ are
contracted into the points $\phi(\Pi_{i,j})=p_{i,j}\in
Z\setminus\widetilde{Q}$. Every fiber of $\widetilde{\phi}$ is of
dimension less or equal to 1 on $V=\p^6\setminus
(\widetilde{Q}\cup(\cup_{i,j}p_{i,j}))$ so that by theorem
\ref{te:Danilov} the restriction of $\widetilde{\phi}$ to
$U=\widetilde{\phi}^{-1}(V)$ is the blowing-up of $\p^6$ along the
smooth variety $Z\setminus
(\widetilde{Q}\cup(\cup_{i,j}p_{i,j}))$. By theorem 1.1 of
\cite{ESB}, we also know that the points $p_{i,j}$ and the surface
$\widetilde{Q}$ are the singular locus of $Z$, so that $Z=X^2$. By
reasoning as in the case discussed in section 4, one immediately
deduces that $X^2$ is the base locus scheme of $\phi^{-1}$ and
that $\phi^{-1}$ is the rational map given by the linear systems
of quartics containing $X^2$, by lemma \ref{equations} (see
\cite{ST2} pg. 207). Moreover we also get that the points
$p_{i,j}$ are isolated triple points for $X^2$, that
$\widetilde{Q}$ is the locus of double points for $X^2$, and that
the planes $\Pi_{i,j}$ are loci of points of multiplicity 2 of the
strict transform of $\Sec(X^1)$, $E_2=\widetilde{ Sec(X^1)}$, i.e.
of the exceptional divisor of $\widetilde{\phi}$. The tangent cone
of $X^2$ at the points $p_{i,j}$ is easily seen to be a cone over
a Segre 3-fold $\p^1\times\p^2\subset\p^5$, see also \cite{ST2}.
 Now, by taking two general quartic hypersurfaces $Y_1$ and $Y_2$
 passing through $X^2$, we have that $Y_1\cap Y_2=X^2\cup (X^2)'$, with $(X^2)'$
 the image of the restriction of $\phi$ to a general $\p^4$ inside the original $\p^6$.
 Since $X^1$ has degree 8, the linear system
 giving the restriction is a linear system of quadrics in $\p^4$ passing through 8 points
 in general linear position, so that $(X^2)'$ is the projection
 of $\nu_2(\p^4)$ from the linear span of 8 general points on it and $\deg(X^2)=\deg((X^2)')=8$.
\end{proof}
 We remark that, see \cite{ST2}, one can also show that $X^2$ is isomorphic
to the projection of $\nu_2(\p^4)$ from the linear span of 8
suitable points on it; it is in fact a degree 8 variety in $\p^6$
having exactly the same singularities as $X^2$. Then the singular
double points of $X^2$ on the quadric surface
$\widetilde{Q}\subset X^2$  correspond to pairs of points in
$\p^4$ imposing only one linear condition to the linear system of
hyperquadrics in $\p^4$
  passing through the 8 points spanning the center of the
  projection
  of $\nu_2(\p^4)$.

\section{Cubic hypersurfaces through a degree 8 3-fold with one
apparent double point}

Once we have done the effort of studying in detail the geometry of
the secant variety of the octic surface in the previous example,
at least  outside the cone over the Segre 3-fold, we can do the
same for another interesting smooth variety $X\subset\p^7$ with
one apparent double point and whose general hyperplane section is
the above mentioned surface.

The construction is the following. Take a line $l\subset \p^7$ and
a skew $\p^5$ inside the fixed $\p^7$. Let $M\subset \p^7$ be the
cone with vertex $l$ above a Segre 3-fold $\p^1\times\p^2$ inside
the fixed $\p^5$. Take the complete intersection of two general
divisors  of type $(1,2)$ on $M$ (see \cite{CMR}). Notice that a
divisor of type $(1,2)$ on $M$ is once again the intersection of
$M$ with a quadric hypersurface containing a $\p^4$ of the ruling.
The resulting variety is a smooth 3-fold $X\subset\p^7$ of degree
8 such that $\Sec(X)=\p^7$ and such that through the general point
of $\p^7$ there passes a unique secant line to $X$. A general
hyperplane section is clearly an octic  surface of the type
discussed in Semple and Tyrrell example of the previous section,
see also \cite{CMR}. The variety $X$ is arithmetically
Cohen-Macaulay and its ideal is generated by 7 degree 2 forms
giving rise to, a rational map $\phi:\p^7\map\p^6$ not defined exactly along $X$.
By restricting to a general hyperplane section,
it immediately follows that the general fiber of $\phi$ is a line, necessarily a secant line to
$X$, i.e. $|H^0(\I_X(2))|$ is a special subhomaloidal system.
\medskip

\begin{Proposition} {\rm (Semple and Tyrrell conic fibration)}
Let notations be as above, let $\phi:\p^7\map\p^6$ be the rational
map defined by the linear system of quadrics vanishing along $X$
and let $\widetilde{\phi}:Bl_X(\p^7)\to\p^6$ be the resolution of
the indeterminacy of $\phi$. Let $\widetilde{M}$ be the strict
transform of $M$ in $Bl_{X}(\p^7)$. Then for every $p\in
Bl_{X}(\p^7)\setminus \widetilde{M}$ the fiber
$\widetilde{\phi}^{-1}(\widetilde{\phi}(p))$ is isomorphic to a
linear space $\p^k_p$ cutting $X$ along a quadric hypersurface
$Q_p=\p^k_p\cap X\subset\p^k_p$. There are three possibilities:
\begin{enumerate}

\item[(i)] $k=1$, $\p^1_p$ is a secant (or tangent) line to $X$
and $Q_p$ is a couple of points;

\item[(ii)] $k=2$, $p$ belongs to one of the $\p^4$'s generated by
a $S(1,2)\subset X$ and $Q_p$ is the conic cut along $S(1,2)$ by
the plane $\p^2_p$;

\item[(iii)] $k=3$, $p$ belongs to one of the 28 linear spaces
$\Pi_{i,j}\simeq\p^3$, $i<j$, $i=1,\ldots, 8$, generated by the 28
quadrics $Q_{i,j}$ contained in $X$ , $\p^3_p=\Pi_{i,j}$ and
$Q_p=Q_{i,j}$.
\end{enumerate}

For any $p\in \widetilde{M}$, the fiber consists of a three
dimensional quadric hypersurface and every $\p^4$ of the ruling
cuts the variety $X$ along a degree 4 cone over an elliptic normal
curve. In particular, $\widetilde{\phi}:Bl_X(\p^7)\to\p^6$ is a
Fano-Mori contraction of fiber type having 28 isolated
3-dimensional fibers isomorphic to $\p^3$, eight 3-dimensional
families of 2-dimensional fibers isomorphic to $\p^2$, each of
them parameterized by a Segre 3-fold, and  a 2-dimensional family
of 3-dimensional fibers isomorphic to quadric hypersurfaces and
parameterized by a smooth quadric surface.\label{Semplep7}
\end{Proposition}

\begin{proof} By cutting with a general hyperplane section
we get a Cremona transformation of the type described in Semple
and Tyrrell example. The analysis of the previous section yields
that each fiber of $\widetilde{\phi}$ is either linear and of
dimension 1, 2 or 3 or a 3 dimensional quadric hypersurfaces.

The locus of the secant lines through every point $p\in\p^7\setminus M$
is a $\p^5_p$ cutting $M$ along a rank 4 quadric $Q_p$. The secant
lines to $X$ passing through $p$ are the secant lines through $p$
of the intersection of $Q_p$ with two divisors of type $(1,2)$ on $M$,
i.e. a 2-dimensional degree 4 (reduced) variety. For general $p$ it will be a
smooth rational normal scroll of degree 4, which is known to be a
surface with one apparent double point. For special $p$ the
variety can degenerate into the union of a plane and a smooth
degree 3 rational normal scroll $S(1,2)$ intersecting along a line,
or into the union of a quadric surface and two skew planes
intersecting the quadric along a line.
There are 28 quadric surface of this kind generating 28 $\p^3$'s as it is seen by cutting with a general
hyperplane and recalling that a general hyperplane section
contains 28 conics (or directly by representing $X$ as the
$\p^1$-bundle $\p(\E)$ over $\p^2$, where $\E$ is the very ample
vector bundle given by the extension
$0\to\O_{\p^2}\to{\E}\to\I_{\{p_1,\ldots,p_8\},{\p^2}}(4)\to 0$,
$ p_1,\ldots,p_8$ in general position, see \cite{CMR}).
Moreover, the surfaces $S(1,2)\subset X$ are
divided into 8 pencils corresponding to the pencils of
lines passing through the 8 points in the above representation of $X$ as a $\p^1$-bundle. The last
claims follow in the same way as in the previous section.
\end{proof}

Now we can produce our last example of Fano-Mori birational
contraction. By abuse of notation we also call $\widetilde{\phi}$
this contraction.

\begin{Theorem}\label{Fanop7} Let $X\subset\p^7$ be as above
and let $Y\subset\p^7$ be a smooth cubic hypersurface through $X$.
Then $\phi_{|Y}:Y\map\p^6$ is birational. The morphism
$\widetilde{\phi}:\widetilde{Y}=Bl_X(Y)\to\p^6$ is  a Fano-Mori
birational contraction. The general fiber of
$\widetilde{\phi}_{|E_2}:E_2\to Z:=\widetilde{\phi}(E_2)$
corresponds to a secant line to $X$ contained in $Y$,
$Z\subset\p^6$ is a 4-fold of degree 13 and $\phi^{-1}:\p^6\map Y$
is given by the linear system of quintics vanishing along $Z$. The
cubic hypersurface $Y$ contains at most a finite number of
$\p^4$'s spanned by the cubic scrolls $S(1,2)$ contained in $X$,
at most a finite number of  $\p^4$'s generating the cone $M$ and
clearly at most a finite number of the $\Pi_{i,j}$. By $\phi$
every $\p^4$ of the ruling of $M$ contained in $Y$ goes into a
line in $Z$ and this line is a locus of singular points for $Z$;
moreover for each point $p$ in the strict transform of this $\p^4$
the fiber $\widetilde{\phi}^{-1}(\widetilde{\phi}(p))$ is
isomorphic to a 3-dimensional quadric in $\p^4$; in particular
these 3-dimensional fibers vary in a one-dimensional family. By
$\phi$ any $\p^4$ spanned by a cubic scroll $S(1,2)$ and contained
in $Y$ goes into a plane contained in $Z$, and this plane is a
locus of singular points for $Z$; moreover for any point $p$ in
the strict transform of this $\p^4$ the fiber
$\widetilde{\phi}^{-1}(\widetilde{\phi}(p))$ is isomorphic to a
plane. Every $\Pi_{i,j}$ contained in $Y$ contracts into a triple
point for $Z$. The remaining positive dimensional fibers of
$\widetilde{\phi}$ come from secant lines to $X$ contained in $Y$
or from quadric surfaces, which are the residual intersection,
outside $X$, of $Y$ with a 3-dimensional quadric in a $\p^4$ of
the ruling of $M$ not contained in $Y$; moreover these quadric
surfaces vary in a two dimensional family parameterized by an open
subset of a quadric surface. Hence $Z$ has a two dimensional locus
of singular points at least, even if $Y$ does not contain any
$\p^4$ of the ruling of $M$ or any $\p^4$ spanned by a cubic
scroll contained in $X$.
\end{Theorem}

\begin{proof}
If $Y$ contains infinitely many's $\p^4$ of the above mentioned
families then it contains a 5 dimensional rational and singular
scroll in $\p^4$'s of degree 3. In fact, in the first case $Y$
would contain $M$, in the second case note that anyone of the
eight families of degree 3 scrolls $S(1,2)$ is contained in such a
scroll. On the other hand, for instance by the
Grothendieck-Lefchetz theorem, it is clear that $Y$ cannot contain
such a scroll, necessarily as a hyperplane section. The assertions
about fibers of $\widetilde{\phi}$ are then easily proved since we
have only to restrict $\widetilde{\phi}:Bl_X(\p^7)\to\p^6$ to
$Bl_X(Y)$ and the fibers of $\widetilde{\phi}$ are completely
described in proposition \ref{Semplep7}. About the singular points
of $Z$ we can argue as in the previous section. About the last
assertions we recall that the intersection of $X$ with any $\p^4$
of the ruling of $M$ is a cone over an elliptic quartic curve
which is given by the intersection of two quadric cones having
$\l$ as vertex. The quartic surface is the base locus of a pencil
of quadrics and, generically, $Y$ cuts any member of the pencil,
out of the base locus, along a residue quadric surface giving rise
to the quoted 2-dimensional fibres of $\widetilde{\phi}$. So we
have only to prove that $Z$ has degree 13 since  it is clearly of
dimension 4. Lemma \ref{equations2} yields that the inverse map is
given by quintic forms through $Z$. Take two general quintics
forms $W_1, W_2$ through $Z$ and let $W_1\cap W_2=Z\cup Z'$. Then
$Z'$ is a 4-dimensional variety which is the image via $\phi$ of a
general linear section of $Y$ with a general $\p^5$, let us say
$L$. Let $Y'=L\cap Y\subset Y$ and let $C=L\cap X$. Then $C$ is a
smooth degree 8 and genus 3 curve in $\p^5$ and $\phi_{|Y'}$
resolves to a birational morphism
$\widetilde{\sigma}:Bl_{C}(Y')\to Z'\subset\p^6$ given by the
restriction to $Y'$ of the linear system of hyperquadrics passing
through $C$.  Let $H'$ be the pull back in $ Bl_{C}(Y')$ of the
hyperplane class of $Y'$ and let $E'$ be the exceptional divisor
so that $|2H'-E'|$ is the linear system giving rise to
$\widetilde{\sigma}$. Then we have that
$\deg(Z')=(2H'-E')^4=16(H')^4-4(2H')\cdot(E')^3+(E')^4=48-8\deg(C)+3\deg(C)+2g(C)-2=48-64+28=12$
and the claim  $\deg(Z)=13$ follows.
\end{proof}

\end{document}